\documentclass[12pt,twoside]{article}
\usepackage{amssymb,amsmath,amsfonts,mathrsfs}
\usepackage{amscd}
\usepackage[mathlines]{lineno}
\usepackage{graphicx}
\usepackage{float}
\usepackage{multirow}
\usepackage{subcaption}
\usepackage{color}
\usepackage{longtable}
\usepackage[numbers,sort&compress]{natbib}
\usepackage[bookmarksnumbered=true,colorlinks,linkcolor=black]{hyperref}
\hypersetup{hypertex=true,colorlinks=true,linkcolor=red,anchorcolor=blue,citecolor=red}
\hypersetup{CJKbookmarks=true}

\newtheorem{thm}{Theorem}[section]
\newtheorem{rmk}{Remark}[section]

\newtheorem{eg}{Example}[section]

\newcounter{saveeqn}%

\makeatletter
\evensidemargin
\oddsidemargin
\makeatother

\title{\Large\bf
Dynamics of polynomial generalized Li\'{e}nard system near the origin and infinity
\thanks{Supported by NSFC \#12101087 and NSFSPC \#2024NSFSC1400.}
}

\author{
{\sc Jun Zhang}
\footnote{Corresponding Author. E-mail address: mathzhangjun@163.com}
\\
$^a${\small College of Mathematics and Physics \& Sichuan Geomath Key Lab}\\
{\small Chengdu University of Technology, Sichuan 610059, P. R. China}
}
\date{}

\begin{document}

\maketitle

\begin{abstract}
We classify all topological phase portraits of the polynomial generalized Li\'{e}nard system,
determined by three arbitrary polynomials,
at the origin and the infinity.
This yields a complete characterization of monodromy at the origin and the infinity.
Moreover,
we obtain a necessary and sufficient condition for the local center via Cherkas' method.
When the origin is the only equilibrium and it is a center,
there are 5 global phase portraits, including two types of global center.
Further, we prove that the global center is not isochronous.

\vskip 0.2cm


\vskip 0.2cm
{\bf AMS (2020) subject classification:}
34C05, 
34C25, 
34D06, 
34D23. 
\end{abstract}

\baselineskip 15pt    
\parskip 10pt         

\thispagestyle{empty}
\setcounter{page}{1}

\section{Introduction}
\setcounter{equation}{0}
\setcounter{lm}{0}
\setcounter{thm}{0}
\setcounter{rmk}{0}
\setcounter{df}{0}
\setcounter{cor}{0}
\setcounter{pro}{0}

The Li\'{e}nard system,
one of the most important mechanical systems,
is of the form
$
\ddot x+f(x) \dot x+g(x)=0,
$
where $g$ represents the restoring force and
$f$ denotes the friction coefficient.
It can be rewritten as a planar differential system
\begin{align}
\dot x=y-F(x),~~~\dot y=-g(x),
\label{equ:Lie}
\end{align}
where $F(x):=\int_0^x f(s)ds$.
Since Li\'{e}nard \cite{Lie28} established a result
on the existence and uniqueness of a limit cycle in 1928,
much attention has been paid to system~\eqref{equ:Lie} on various topics:
limit cycle,
invariant algebraic curve,
integrability,
center,
isochronicity,
global center,
phase portrait,
bifurcation,
etc.
Moreover,
different kinds of generalization of Li\'{e}nard system are also proposed and studied
(see \cite{CX23,GWZ,GT98,GT99,Gine,VZ,XZ}).
In this paper,
we focus on the polynomial generalization
\begin{align}
\dot x=\varphi(y)-F(x),~~~\dot y=-g(x),
\label{GL:pfg}
\end{align}
where
$\varphi(y):=a_p y^p+\cdots+a_\ell y^\ell$,
$F(x):=b_q x^q+\cdots+b_m x^m$,
$g(x):=c_r x^r+\cdots+c_n x^n$,
$a_pb_qc_ra_\ell b_mc_n\ne 0$,
$p\le \ell$, $q\le m$, $r\le n$,
and $p,q,r\ge 1$.
Note that system~\eqref{GL:pfg} becomes system~\eqref{equ:Lie} when $\ell=1$.
Compared to system~\eqref{equ:Lie},
system~\eqref{GL:pfg} is determined by three arbitrary polynomials
and the origin can be a degenerate equilibrium with zero linear part when $p,q,r\ge 2$.

Gasull and Torregrosa \cite{GT98} proved that
the equilibrium $O$ of system~\eqref{GL:pfg} is monodromic
only when $p$ and $r$ are both odd and $q(p+1)-p(r+1)\ge 0$,
and solved the monodromy problem completely in the subcase $q(p+1)-p(r+1)>0$
by introducing Lyapunov polar coordinates.
As indicated in \cite[Remark~2.3]{GT98},
the subcase $q(p+1)-p(r+1)=0$ is much more delicate,
so they studied the situation $p=1$ only.
Moreover,
equivalent necessary and sufficient conditions for system~\eqref{GL:pfg}
to have a center at the origin are obtained.
Further,
they applied their results to solve the center problem
for several families of polynomial differential systems,
which especially provides a new and unified way to
characterize all centers for quadratic systems.
Later,
under the above monodromic condition,
using the multiplicity of a polynomial map generated by $F$ and the primitive of $g$,
Gasull and Torregrosa \cite{GT99} provided an upper estimate on
the number of small-amplitude limit cycles
bifurcating from the weak focus $O$,
Similar result can also be found in \cite{Han}.
However,
monodromy condition and center condition in the subcase $q(p+1)-p(r+1)=0$
are left to be solved,
which would extend the range of the application of Gasull and Torregrosa's idea
to center problem.
Meanwhile,
a complete classification of topological phase portraits
at the origin remains to be answered.

On the other hand,
investigating behaviour of orbits for polynomial differential systems near infinity
is meaningful in several aspects.
Clearly, it is a prerequisite step to study the global phase portraits.
One can also benefit from it to construct outer boundary curve
when detecting limit cycles via Poincar\'{e}-Bendixson Theorem.
It helps us search algebraic invariant curves
and also detect centers having infinity in the boundary of their period annulus.
Moreover,
it yields some necessary conditions for isochronicity of a center.
Dumortier and Herssens \cite{DH99} gave a complete classification of phase portraits of
the polynomial Li\'{e}nard system
\begin{align}
\dot x=y,~~~\dot y=-\tilde{g}(x)-\tilde{f}(x)y
\label{equ:DH}
\end{align}
near infinity,
except for the center-focus identification.
Note that system \eqref{equ:DH} can be changed into the form of system~\eqref{equ:Lie}
by the transformation $y\to y+\int_0^x \tilde{f}(s)ds$.
Similarly to polynomial Li\'{e}nard system~\eqref{equ:DH},
Chen, Zhang and Zhang \cite{CZZ}
characterized dynamics of the polynomial Rayleigh-Duffing system
\begin{align}
\dot x=y,~~~\dot y=-\tilde{g}(x)-\tilde{f}(y)y
\label{equ:RD}
\end{align}
near infinity,
obtained a necessary and sufficient condition for an equilibrium to be a center,
and classified all global phase portraits of the Rayleigh-Duffing system
on the Poincar\'{e} disc when there is only one equilibrium and it is a center.
Schlomiuk and Vulpe \cite{SN05}
investigated the geometry of quadratic differential systems
in a neighborhood of the infinity, see more related results in \cite{ALSV}.
Recently,
given a planar polynomial differential system with a fixed Newton polytope,
Dalbelo, Oliveira and Perez \cite{DOP} used monomials associated to
the upper boundary of the Newton polytope
to determine the topological phase portrait near infinity.
However,
dynamics of system~\eqref{GL:pfg} at infinity have not been studied yet.

Once the center condition at the origin and the monodromy condition at infinity are obtained,
one can further determine whether the local center is a global one.
That is an open problem proposed by Conti \cite[pp.228-230]{Conti98}:
{\it Identify all polynomial systems having a $($isochronous$)$ global center.}
Llibre and Valls \cite{LV22} provided a description on
the nondegenerate global center of system~\eqref{equ:DH}.
Using phase portraits of system~\eqref{equ:DH} near infinity obtained in \cite{DH99},
Chen, Li and Zhang \cite{CLZ} gave a necessary and sufficient condition for
(nondegenerate and nilpotent) global center.
It is proved in \cite{CZZ} that system~\eqref{equ:RD} has no global centers.
Global centers for some generalized Kukles polynomial systems
and other systems can be found in recent works \cite{CFZ, CL24} and references therein.
For isochronicity of global center,
it is proved that
nonlinear polynomial rigid system (\cite{Conti94}),
potential system (\cite[Theorem~4.1]{CJ89}),
separable Hamiltonian system (\cite{CGM})
and polynomial system with linear plus homogeneous nonlinearities (\cite{LSZ})
have no isochronous global center.
Authors in \cite{LSZ} obtained all cubic polynomial differential systems
having an isochronous global center.
For systems of higher degree,
Authors in \cite{CMV99} characterized
isochronous global center for Hamiltonian system of the form
$H(x,y)=A(x)+B(x)y+C(x)y^2$ with polynomials $A$, $B$ and $C$.
However,
global center and its isochronicity
for generalized Li\'{e}nard system~\eqref{GL:pfg} have not been investigated
and it would be more complicated compared to Li\'{e}nard system \eqref{equ:DH}
since there is an extra arbitrary polynomial $\varphi$.

In this paper
we investigate dynamics of system~\eqref{GL:pfg}.
In section 2,
using quasi-homogeneous blow-up and
computing center manifolds for semi-hyperbolic equilibria (of the desingularized vector field)
with the aid of Newton polygon,
we prove that there are 9 different topological phase portraits of
system~\eqref{GL:pfg} near the origin
and everything is determined by $a_p,b_q, c_r$ and parities of $p,q$ and $r$,
except for the center-focus identification.
This also suggests a complete characterization of monodromy at the origin,
which is partially solved in \cite{GT98}.
Moreover,
a necessary and sufficient condition of the local center is obtained by the Cherkas' method.
Section 3 is devoted to dynamics of system~\eqref{GL:pfg} with $\ell\ge 2$ near infinity
since the case $\ell=1$ has already been studied in \cite{DH99}.
By successive blow-ups and
checking how the blow-up affects the Newton polygon,
we prove that
there are 26 different topological phase portraits and all information is determined by
$a_{\ell}, b_m$ and $c_n$ and parities of $\ell, m$ and $n$,
except for the center-focus problem at infinity.
In section 4
we analyze the boundary of the period annulus and use Poincar\'{e}-Hopf Index Theorem to
classify all global phase portraits of system~\eqref{GL:pfg} with $\ell\ge 2$
on the Poincar\'{e} disc when the origin is the only equilibrium and it is a center.
Moreover, a necessary and sufficient condition of the global center is obtained.
Further, the global center is proved to be non-isochronous
by investigating the order of the meromorphic vector field
obtained by desingularization near infinity.
Thus Conti's open problem is solved for generalized Li\'{e}nard system~\eqref{GL:pfg}.

\section{Phase portraits at the origin}
\setcounter{equation}{0}
\setcounter{lm}{0}
\setcounter{thm}{0}
\setcounter{rmk}{0}
\setcounter{df}{0}
\setcounter{cor}{0}
\setcounter{pro}{0}

In this section
we classify phase portraits of system~\eqref{GL:pfg} near the origin.
Due to the time-rescaling $t\to -t/|a_p|$,
we can assume without loss of generality that $a_p=-1$.
Although we consider topological phase portraits,
we still distinguish nodes and foci in the following
since a node has characteristic orbits and it is not monodromic.

\begin{figure}[H]
\centering
\subcaptionbox{%
     }{\includegraphics[height=1in]{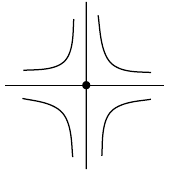}}~
\subcaptionbox{%
     }{\includegraphics[height=1in]{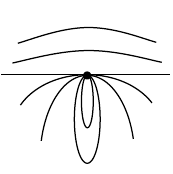}}~
\subcaptionbox{%
     }{\includegraphics[height=1in]{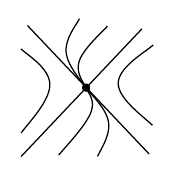}}~
\subcaptionbox{%
     }{\includegraphics[height=1in]{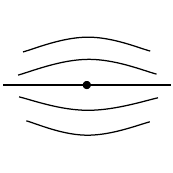}}
     \\
\subcaptionbox{%
     }{\includegraphics[height=1in]{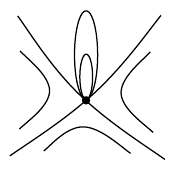}}~
\subcaptionbox{%
     }{\includegraphics[height=1in]{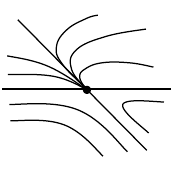}}~
\subcaptionbox{%
     }{\includegraphics[height=1in]{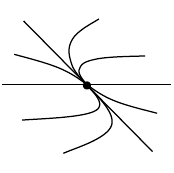}}~
\subcaptionbox{%
     }{\includegraphics[height=1in]{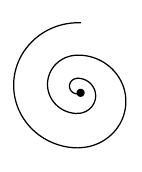}}~
\subcaptionbox{%
     }{\includegraphics[height=1in]{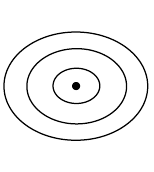}}
\caption{Topological phase portrait of system~\eqref{GL:pfg} at the origin.}
\label{fig:O-phase}
\end{figure}

\begin{thm}
System~\eqref{GL:pfg} with $a_p=-1$ has 9 different topological phase portraits at the origin,
see Fig.~\ref{fig:O-phase}.
More concretely,
the relation between the phase portraits and parameters is listed
in Table~\ref{tab:O}, in which $\hat{c}:=q\big|\frac{b_q}{p+1}\big|^{\frac{p+1}{p}}$.
\label{thm:O}
\end{thm}

\begin{center}
\begin{table}[H]
\small
\begin{tabular}{c|c|r|r|r|c|l}
\hline
Fig.\ref{fig:O-phase}(a)
&                & odd $p$  &          & odd $r$  &         & $c_r>0$ \\
\hline
Fig.\ref{fig:O-phase}(h), (i)
&\multirow{3}{*}{$p(r+1)<q(p+1)$}
                         & odd $p$  &          & odd $r$  &         & $c_r<0$ \\
Fig.\ref{fig:O-phase}(d)
&                & odd $p$  &          & even $r$ &         &         \\
Fig.\ref{fig:O-phase}(d)
&                & even $p$ &          &          &         &         \\
\hline
Fig.\ref{fig:O-phase}(g)
&\multirow{7}{*}{$p(r+1)=q(p+1)$}
                         & odd $p$  & odd $q$  & odd $r$  &         & $c_r\in[-\hat{c},0)$ \\
Fig.\ref{fig:O-phase}(h)(i)
&                & odd $p$  &          & odd $r$  &         & $c_r<-\hat{c}$      \\
Fig.\ref{fig:O-phase}(b)
&                & odd $p$  & even $q$ & odd $r$  &         & $c_r\in[-\hat{c},0)$ \\
Fig.\ref{fig:O-phase}(d)
&                & even $p$ & even $q$ &          & $b_q>0$ &                \\
Fig.\ref{fig:O-phase}(d)
&                & even $p$ & even $q$ &          & $b_q<0$ & $|c_r|>\hat{c}$\\
Fig.\ref{fig:O-phase}(e)
&                & even $p$ & even $q$ & odd $r$  & $b_q<0$ & $|c_r|\le \hat{c}$\\
Fig.\ref{fig:O-phase}(c)
&                & even $p$ & even $q$ & even $r$ & $b_q<0$ & $|c_r|\le \hat{c}$\\
\hline
Fig.\ref{fig:O-phase}(g)
&\multirow{6}{*}{$p(r+1)>q(p+1)$}
                         & odd $p$  & odd $q$  & odd $r$  &         & $c_r<0$ \\
Fig.\ref{fig:O-phase}(b)
&                & odd $p$  & even $q$ & odd $r$  &         & $c_r<0$ \\
Fig.\ref{fig:O-phase}(f)
&                & odd $p$  &          & even $r$ &         &         \\
Fig.\ref{fig:O-phase}(d)
&                & even $p$ & even $q$ &          & $b_q>0$ &         \\
Fig.\ref{fig:O-phase}(e)
&                & even $p$ & even $q$ & odd $r$  & $b_q<0$ &         \\
Fig.\ref{fig:O-phase}(c)
&                & even $p$ & even $q$ & even $r$ & $b_q<0$ &         \\
Fig.\ref{fig:O-phase}(f)
&                & even $p$ & odd $q$  &          &         &         \\
\hline
\end{tabular}
\caption{Parameter conditions for different topological phase portraits.}
\label{tab:O}
\end{table}
\end{center}

\vspace{-50pt}

The proof of Theorem~\ref{thm:O} will be given latter because of its length.
Noticing that center and focus are not identified when $p(r+1)\le q(p+1)$ in Theorem~\ref{thm:O},
so we characterize the monodromy and center in the following theorem.

\begin{thm}
System~\eqref{GL:pfg} with $a_p=-1$ is monodromic at $O$ if and only if either
\\
{\bf (M1)} $p(r+1)<q(p+1)$, both $p$ and $r$ are odd, and $c_r<0$; or
\\
{\bf (M2)} $p(r+1)=q(p+1)$, both $p$ and $r$ are odd, and $c_r<-\hat{c}$.
\\
Moreover,
the monodromic equilibrium $O$ is a center
if and only if the system $F(x)=F(z)$ and $G(x)=G(z)$
has a unique solution $z(x)$ satisfying $z(0)=0$ and $z'(0)<0$,
where $G(x):=\int_0^x g(s) ds$.
\label{thm:center}
\end{thm}

Necessary and sufficient condition for monodromy is obtained from Theorem~\ref{thm:O} directly.
When the equilibrium $O$ is monodromic,
necessary and sufficient condition for local center can be obtained by Cherkas' method,
the same as the proof of \cite[Theorem~2.6]{GT98} and \cite[Theorem~6]{Ch99},
in both which there are some other equivalent necessary and sufficient conditions.

\begin{figure}[!h]
\centering
\subcaptionbox{%
     }{\includegraphics[height=1in]{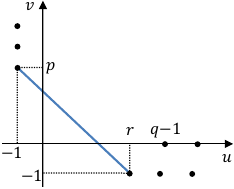}}~~~~~~
\subcaptionbox{%
     }{\includegraphics[height=1in]{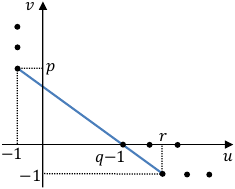}}~~~~~~
\subcaptionbox{%
     }{\includegraphics[height=1in]{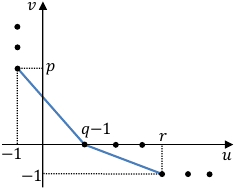}}
\caption{Three cases of Newton polygons of system~\eqref{GL:pfg}.}
\label{fig:NP}
\end{figure}

{\bf Proof of Theorem~\ref{thm:O}.}
In order to obtain the topological phase portraits
of the equilibrium $O$ of system~\eqref{GL:pfg},
we use quasi-homogeneous blowing-up associated with Newton polygon
to desingularize the equilibrium $O$ (see \cite[Chapter~3]{DLA}).
We see from the expansion of system~\eqref{GL:pfg} that
all support points lie on either the line $u=-1$,
or the line $v=0$ or the line $v=-1$.
So the Newton polygon of system~\eqref{GL:pfg} has three cases:
{\bf (I)} only one edge, linking vertices $(-1,p)$ and $(r,-1)$, and
there are no support points on the edge except for vertices,
{\bf (II)} only one edge, linking vertices $(-1,p)$ and $(r,-1)$,
and the point $(q-1,0)$ is the only support point on the edge except for vertices, and
{\bf (III)} exactly two edges,
linking vertices $(-1,p)$, $(q-1,0)$ and $(r,-1)$,
as illustrated in Fig.~\ref{fig:NP} (a)-(c), respectively.

\begin{figure}[h]
\centering
\subcaptionbox{%
     }{\includegraphics[height=1in]{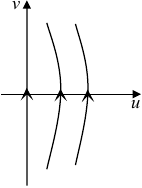}}~~~~~~
\subcaptionbox{%
     }{\includegraphics[height=1in]{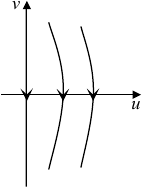}}
\caption{Phase portraits of systems~\eqref{equ:OI+} and \eqref{equ:OI-}
along the $v$-axis.}
\label{fig:OI}
\end{figure}

In case {\bf(I)}, we have $p(r+1)<q(p+1)$,
and the only edge lies on the line $\alpha u+\beta v=\delta_1:=pr-1$,
where $\alpha:=p+1$ and $\beta:=r+1$.
As done in \cite[Chapter~3]{DLA},
we rewrite system~\eqref{GL:pfg}
in quasi-homogeneous components of type $(\alpha,\beta)$ as
\begin{equation*}
\dot x=P^{(1)}_{\delta_1}(x,y)+\cdots,~~~
\dot y=Q^{(1)}_{\delta_1}(x,y)+\cdots,
\end{equation*}
where $P^{(1)}_{\delta_1}(x,y):=-y^p$, $Q^{(1)}_{\delta_1}(x,y)=-c_rx^r$,
and dots represent those terms of quasi-homogeneous degree bigger than $\delta_1$.
Blowing up the equilibrium $O$ in the positive $x$-direction
by the transformation $x=u^\alpha$ and $y=u^\beta v$, we obtain
\begin{equation}
\dot u=u\{P^{(1)}_{\delta_1}(1,v)+O(u)\},~~~
\dot v=G_0(v)+O(u),
\label{equ:OI+}
\end{equation}
where a time-rescaling is performed,
$P^{(1)}_{\delta_1}(1,v)=-v^p$ and
$$
G_0(v):=\alpha Q^{(1)}_{\delta_1}(1,v)-\beta vP^{(1)}_{\delta_1}(1,v)
=-(p+1) c_r+(r+1) v^{p+1}.
$$
If system~\eqref{equ:OI+} has an equilibrium $(0,v_*)$ on the $v$-axis,
then it is a hyperbolic one with the Jacobian matrix
\begin{equation*}
\left(
\begin{array}{ccc}
-v_*^p  &  0
\\
\star   &  G'_0(v_*)
\end{array}
\right).
\end{equation*}
On the other hand,
blowing up the equilibrium $O$ of system~\eqref{GL:pfg}
in the negative $x$-direction
by the transformation $x=-u^\alpha$ and $y=u^\beta v$, we obtain
\begin{equation}
\dot u=-u\{P^{(1)}_{\delta_1}(-1,v)+O(u)\},~~~
\dot v=\widetilde{G}_0(v)+O(u),
\label{equ:OI-}
\end{equation}
where a time-rescaling is performed,
$P^{(1)}_{\delta_1}(-1,v)=-v^p$ and
$$
\widetilde{G}_0(v):=\alpha Q^{(1)}_{\delta_1}(-1,v)+\beta vP^{(1)}_{\delta_1}(-1,v)
=(-1)^{r+1} (p+1) c_r -(r+1) v^{p+1}.
$$
If system~\eqref{equ:OI-} has an equilibrium $(0,\tilde{v}_*)$ on the $v$-axis,
then it is a hyperbolic one with the Jacobian matrix
\begin{equation*}
\left(
\begin{array}{ccc}
\tilde{v}_*^p  &  0
\\
\star          &  \widetilde{G}'_0(\tilde{v}_*)
\end{array}
\right)
\end{equation*}
Since equilibria of systems~\eqref{equ:OI+} and \eqref{equ:OI-} on the $v$-axis
depend on the parities of $p$ and $r$ and the sign of $c_r$,
we consider the following 5 subcases:
{\bf(I1)} $p$ is even,
{\bf(I2)} $p$ is odd, $r$ is odd,  $c_r>0$,
{\bf(I3)} $p$ is odd, $r$ is odd and $c_r<0$,
{\bf(I4)} $p$ is odd, $r$ is even and $c_r>0$, and
{\bf(I5)} $p$ is odd, $r$ is even and $c_r<0$.
We only give the proof in subcase {\bf(I3)} since the proofs in other subcases are similar. In {\bf(I3)},
neither system~\eqref{equ:OI+} nor system~\eqref{equ:OI-}
has an equilibrium on the $v$-axis.
So phase portraits of systems~\eqref{equ:OI+} and \eqref{equ:OI-}
along the $v$-axis in the half-plane $u\ge 0$ are given by
Fig.~\ref{fig:OI} (a) and (b), respectively.
After blowing down,
we find that the equilibrium $O$ of system~\eqref{GL:pfg} is monodromic,
having phase portrait Fig.~\ref{fig:O-phase} (h) or (i).

\begin{figure}[h]
\centering
\subcaptionbox{%
     }{\includegraphics[height=1in]{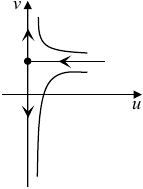}}~~
\subcaptionbox{%
     }{\includegraphics[height=1in]{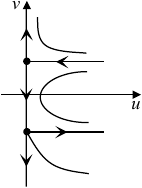}}~~
\subcaptionbox{%
     }{\includegraphics[height=1in]{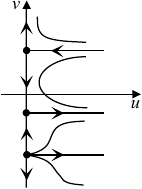}}~~
\subcaptionbox{%
     }{\includegraphics[height=1in]{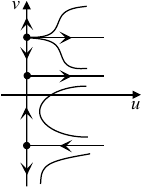}}~~
\subcaptionbox{%
     }{\includegraphics[height=1in]{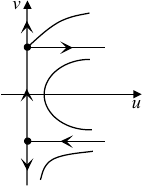}}~~
\subcaptionbox{%
     }{\includegraphics[height=1in]{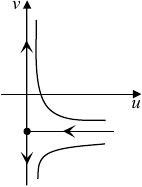}}
\caption{Phase portraits of system~\eqref{equ:OII4+} along the $v$-axis.}
\label{fig:OII4+}
\end{figure}

In case {\bf(II)},
we have $p(r+1)=q(p+1)$, and
the only edge lies on the line $pu+qv=\delta_2:=p(q-1)$.
As done in \cite[Chapter~3]{DLA},
we rewrite system~\eqref{GL:pfg} in quasi-homogeneous components of type $(p,q)$ as
\begin{eqnarray*}
\dot x=P^{(2)}_{\delta_2}(x,y)+\cdots,~~~
\dot y=Q^{(2)}_{\delta_2}(x,y)+\cdots,
\end{eqnarray*}
where $P^{(2)}_{\delta_2}(x,y):=-y^p+b_qx^q$, $Q^{(2)}_{\delta_2}(x,y):=-c_rx^r$
and dots represent terms of quasi-homogeneous degree bigger than $\delta_2$.
Then blowing up the equilibrium $O$ of system~\eqref{GL:pfg}
in the positive $x$-direction by the transformation $x=u^p$ and $y=u^qv$,
we obtain
\begin{eqnarray}
\dot u=u\{P^{(2)}_{\delta_2}(1,v)+O(u)\},~~~
\dot v=H_0(v)+O(u),
\label{equ:OII4+}
\end{eqnarray}
where a time-rescaling is performed,
$P^{(2)}_{\delta_2}(1,v)=-b_q-v^p$ and
$$
H_0(v):=pQ^{(2)}_{{\delta_2}}(1,v)-qvP^{(2)}_{{\delta_2}}(1,v)
=-pc_r+qb_q v+q v^{p+1}.
$$
On the other hand,
blowing up the equilibrium $O$ of system~\eqref{GL:pfg}
in the negative $x$-direction by the transformation $x=-u^p$ and $y=u^qv$,
we obtain
\begin{eqnarray}
\dot u=-u\{P^{(2)}_{\delta_2}(-1,v)+O(u)\},~~~
\dot v=\widetilde{H}_0(v)+O(u),
\label{equ:OII4-}
\end{eqnarray}
where a time-rescaling is performed,
$P^{(2)}_{\delta_2}(-1,v)=(-1)^{q+1}b_q-v^p$ and
$$
\widetilde{H}_0(v)
:=pQ^{(2)}_{{\delta_2}}(-1,v)+qvP^{(2)}_{{\delta_2}}(-1,v)
=(-1)^{r+1}pc_r+(-1)^{q+1}qb_qv-qv^{p+1}.
$$
Equilibria of systems~\eqref{equ:OII4+} and \eqref{equ:OII4-} on the $v$-axis
are determined by real zeros of polynomials $H_0$ and $\widetilde{H}_0$,
which depend on the parities of $p,q$ and $r$ and also the sign of $b_q$.
Note that the equality $p(r+1)=q(p+1)$ implies that
$r$ is odd when $p$ is odd and $q$ is even when $p$ is even.
Then there are 6 subcases:
\begin{center}
\begin{tabular}{ll}
{\bf(II1)} odd $p,q$ and $r$,               & {\bf(II2)} odd $p$, even $q$ and odd $r$,  \\
{\bf(II3)} even $p,q$, odd $r$ and $b_q>0$, & {\bf(II4)} even $p,q$, odd $r$ and $b_q<0$, \\
{\bf(II5)} even $p,q,r$, and $b_q>0$,       & {\bf(II6)} even $p,q,r$, and $b_q<0$.       \\
\end{tabular}
\end{center}
We only give the proof in subcase {\bf(II4)}
since proofs in other subcases are similar.
In {\bf(II4)},
the derivative $H'_0(v)$ has two real zeros
$$
v_\pm:=\pm\left(\frac{-b_q}{p+1}\right)^{1/p}.
$$
It follows that $H_0(v)$ has either
{\bf (Z1)} only one real zero ($>v_+$) if $c_r>\hat{c}$, or
{\bf (Z2)} one double zero $v_-$ and one positive zero $(>v_+)$ if $c_r=\hat{c}$, or
{\bf (Z3)} two negative real zeros and one positive one if $c_r\in(0,\hat{c})$, or
{\bf (Z4)} one negative real zero and two positive ones if $c_r\in(-\hat{c},0)$, or
{\bf (Z5)} one negative zero and one double zero $v_+$ if $c_r=-\hat{c}$, or
{\bf (Z6)} only one real zero ($<v_-$) if $c_r<-\hat{c}$.
If $(0,v_*)$ is an equilibrium of system~\eqref{equ:OII4+},
then the Jacobian matrix at this equilibrium is given by
\begin{eqnarray*}
\left(
\begin{array}{cc}
P^{(2)}_{\delta_2}(1,v_*) & 0
\\
\star               & H_0'(v_*)
\end{array}
\right)
=
\left(
\begin{array}{cc}
-\frac{pc_r}{qv_*} & 0
\\
\star             & H_0'(v_*)
\end{array}
\right).
\end{eqnarray*}
Hence phase portraits of system~\eqref{equ:OII4+} along the $v$-axis
in the half-plane $u\ge 0$
in situations {\bf(Z1)}-{\bf(Z6)} are given by
Fig.~\eqref{fig:OII4+} (a)-(f), respectively.

\begin{figure}[H]
\centering
\subcaptionbox{%
     }{\includegraphics[height=1in]{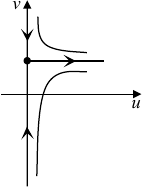}}~~
\subcaptionbox{%
     }{\includegraphics[height=1in]{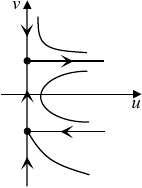}}~~
\subcaptionbox{%
     }{\includegraphics[height=1in]{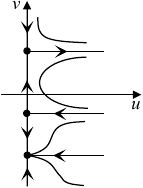}}~~
\subcaptionbox{%
     }{\includegraphics[height=1in]{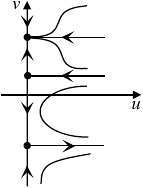}}~~
\subcaptionbox{%
     }{\includegraphics[height=1in]{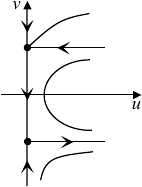}}~~
\subcaptionbox{%
     }{\includegraphics[height=1in]{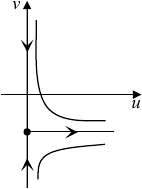}}
\caption{Phase portraits of system~\eqref{equ:OII4-} along the $v$-axis.}
\label{fig:OII4-}
\end{figure}

\noindent
On the other hand,
the derivate $\widetilde{H}'_0(v)$ also has the two real zeros $v_\pm$.
So $\widetilde{H}_0(v)$ also satisfies one of situations {\bf (Z1)}-{\bf (Z6)}
as $H_0(v)$.
If $(0,\tilde{v}_*)$ is an equilibrium of system~\eqref{equ:OII4-},
the Jacobian matrix at this equilibrium is given by
\begin{eqnarray*}
\left(
\begin{array}{cccc}
-P^{(2)}_{\delta_2}(-1,\tilde{v}_*) & 0
\\
\star                         & \widetilde{H}'_0(\tilde{v}_*)
\end{array}
\right)
=
\left(
\begin{array}{cccc}
\frac{pc_r}{q\widetilde{v}_*} & 0
\\
\star                          & \widetilde{H}'_0(\tilde{v}_*)
\end{array}
\right).
\end{eqnarray*}
Thus phase portraits of system~\eqref{equ:OII4-} along the $v$-axis
in the half-plane $u\ge 0$
in situations {\bf(Z1)}-{\bf(Z6)} are given by
Fig.~\ref{fig:OII4-} (a)-(f), respectively.

\noindent
Combining those phase portraits of systems~\eqref{equ:OII4+} and \eqref{equ:OII4-},
we obtain the phase portraits Fig.~\ref{fig:O-phase} (d), (e), (e), (e), (e) and (d)
of system~\eqref{GL:pfg} at the origin in situation {\bf(Z1)}-{\bf(Z6)}, respectively.

In case {\bf(III)},
we have $p(r+1)>q(p+1)$ and
the first edge lies on the line $pu+qv=\delta_2$.
Similarly to case {\bf(I)},
we rewrite system~\eqref{GL:pfg} in quasi-homogeneous components of type $(p,q)$ as
\begin{equation*}
\left\{
\begin{array}{lll}
\dot x=P^{(3)}_{\delta_2}(x,y)+\cdots+P^{(3)}_{\delta_*-1}(x,y)+P^{(3)}_{\delta_*}(x,y)+\cdots,
\\
\dot y=Q^{(3)}_{\delta_2}(x,y)+\cdots+Q^{(3)}_{\delta_*-1}(x,y)+Q^{(3)}_{\delta_*}(x,y)+\cdots
\end{array}
\right.
\end{equation*}
such that
$P^{(3)}_k(r^px,r^qy)=r^{p+k}P^{(3)}_k(x,y)$ and
$Q^{(3)}_k(r^px,r^qy)=r^{q+k}Q^{(3)}_k(x,y)$,
where $P^{(3)}_{\delta_2}(x,y):=-y^p-b_qx^q$, $Q^{(3)}_{\delta_2}(x,y)=0$, $\delta_*:=pr-q$,
and dots represent those terms of quasi-homogeneous degree bigger than $\delta_*$.
Compared to case {\bf(I)},
the desingularized system in case {\bf(III)} may have semi-hyperbolic equilibria,
so more terms are involved in the above expression to
determine their qualitative properties.
Note that
all support points of $\dot y$ in system~\eqref{GL:pfg} lie on the line $v=-1$,
on which $(r,-1)$ is the only support point $(i,j)$ such that
$pi+qj$ takes its minimum $\delta_*$.
Then
$$
Q^{(3)}_{\delta_2+1}(x,y)\cdots=Q^{(3)}_{\delta_*-1}(x,y)=0~~~\mbox{and}~~~
Q^{(3)}_{\delta_*}(x,y)=-c_rx^r.
$$
Blowing up the equilibrium $O$ in the positive $x$-direction
by the transformation $x=u^p$ and $y=u^qv$, we obtain
\begin{equation}
\left\{
\begin{array}{llll}
\dot u={\cal U}(u,v)
:=u\{W(u,v)+P^{(3)}_{\delta_*}(1,v)u^{\delta_*-\delta_2}+O(u^{\delta_*-\delta_2+1})\},
\\
\dot v={\cal V}(u,v)
:=-qvW(u,v)+u^{\delta_*-\delta_2}\{-pc_r-qvP^{(3)}_{\delta_*}(1,v)\}
+O(u^{\delta_*-\delta_2+1}),
\end{array}
\right.
\label{equ:OIII+}
\end{equation}
where a time-rescaling is performed and
$$
W(u,v):=P^{(3)}_{\delta_2}(1,v)+uP^{(3)}_{\delta_2+1}(1,v)+\cdots
+u^{\delta_*-\delta_2-1}P^{(3)}_{\delta_*-1}(1,v).
$$
Clearly,
equilibria of system~\eqref{equ:OIII+} on the $v$-axis are determined by the equation
$$
-qvW(0,v)=-qvP^{(3)}_{\delta_2}(1,v)=qv(v^p+b_q)=0.
$$
Thus system~\eqref{equ:OIII+} has the equilibrium $(0,0)$ and
at most two other equilibria on the $v$-axis.
The Jacobian matrices of system~\eqref{equ:OIII+} at equilibria $(0,0)$ and $(0,v_*)$
(if exists) are given by
\begin{equation*}
\left(
\begin{array}{cc}
-b_q   & 0
\\
\star  & qb_q
\end{array}
\right)
~~~\mbox{and}~~~
\left(
\begin{array}{cc}
0     & 0
\\
\star & pqv_*^p
\end{array}
\right),
\end{equation*}
respectively.
Clearly, the equilibrium $(0,0)$ is a hyperbolic saddle,
and the equilibrium $(0,v_*)$ (if exists) is semi-hyperbolic.
In order to determine its qualitative property,
by Theorem~7.1 in \cite[p.114]{ZZF},
we solve from the equation ${\cal V}(u,v)=0$ near the semi-hyperbolic equilibrium $(0,v_*)$
that $v=\Lambda(u)=v_*+O(u)$ and
$$
W(u,\Lambda(u))=-\frac{pc_r}{qv_*}u^{\delta_*-\delta_2}-P^{(3)}_{\delta_*}(1,v)u^{\delta_*-\delta_2}
+O(u^{\delta_*-\delta_2+1}).
$$
Substituting it in ${\cal U}(u,v)$, we obtain that
$$
{\cal U}(u,\Lambda(u))=-\frac{pc_r}{qv_*}u^{\delta_*-\delta_2+1}+O(u^{\delta_*-\delta_2+2}).
$$
Then the qualitative property of the semi-hyperbolic equilibrium
$(0,v_*)$ (if exists) is determined by the signs of $pqv_*^p$ and $-pc_r/(qv_*)$ and
the parity of $\delta_*-\delta_2+1$.

On the other hand,
blowing up the equilibrium $O$ in the negative $x$-direction
by the transformation $x=-u^p$ and $y=u^qv$, we obtain
\begin{equation}
\left\{
\begin{array}{llll}
\dot u=\widetilde{{\cal U}}(u,v)
:=-u\{\widetilde{W}(u,v)+P^{(3)}_{\delta_*}(-1,v)u^{\delta_*-\delta_2}
+O(u^{\delta_*-\delta_2+1})\},
\\
\dot v=\widetilde{{\cal V}}(u,v)
:=qv\widetilde{W}(u,v)+u^{\delta_*-\delta_2}\{(-1)^{r+1}pc_r+qvP^{(3)}_{\delta_*}(-1,v)\}
+O(u^{\delta_*-\delta_2+1}),
\end{array}
\right.
\label{equ:OIII-}
\end{equation}
where a time-rescaling is performed and
$$
\widetilde{W}(u,v):=P^{(3)}_{\delta_2}(-1,v)+uP^{(3)}_{\delta_2+1}(-1,v)+\cdots
+u^{\delta_*-\delta_2-1}P^{(3)}_{\delta_*-1}(-1,v).
$$
Thus system~\eqref{equ:OIII-} has the equilibrium $(0,0)$ and
at most two other equilibria on the $v$-axis.
The Jacobian matrices of system~\eqref{equ:OIII-} at equilibria $(0,0)$ and
$(0,\tilde{v}_*)$ (if exists) are given by
\begin{equation*}
\left(
\begin{array}{cc}
(-1)^qb_q   & 0
\\
\star       & (-1)^{q+1}qb_q
\end{array}
\right)
~~~\mbox{and}~~~
\left(
\begin{array}{cc}
0     & 0
\\
\star & -pq\tilde{v}_*^p
\end{array}
\right),
\end{equation*}
respectively.
Clearly, the equilibrium $(0,0)$ is a hyperbolic saddle,
and the equilibrium $(0,\tilde{v}_*)$ (if exists) is semi-hyperbolic.
In order to determine its qualitative proper,
similarly to the above,
we solve from the equation $\widetilde{{\cal V}}(u,v)=0$
near the semi-hyperbolic equilibrium $(0,\tilde{v}_*)$ that
$v=\tilde\Lambda(u)=\tilde{v}_*+O(u)$ and
$$
\widetilde{W}(u,\tilde\Lambda(u))
=(-1)^{r}\frac{pc_r}{q\tilde\Lambda(u)}u^{\delta_*-\delta_2}
-P^{(3)}_{\delta_*}(1,v)u^{\delta_*-\delta_2}
+O(u^{\delta_*-\delta_2+1}).
$$
Substituting it in $\widetilde{{\cal U}}(u,v)$, we obtain that
$$
\widetilde{{\cal U}}(u,\tilde\Lambda(u))
=(-1)^{r+1}\frac{ p c_r}{q\tilde{v}_*}u^{\delta_*-\delta_2+1}+O(u^{\delta_*-\delta_2+2}).
$$
Then the qualitative property of the semi-hyperbolic equilibrium
$(0,\tilde{v}_*)$ is determined by the signs of
$-pq\tilde{v}_*^p$ and $(-1)^{r+1}pc_r/(q\tilde{v}_*)$ and
the parity of $\delta_*-\delta_2+1$.

\begin{figure}[H]
\centering
\subcaptionbox{%
     }{\includegraphics[height=1in]{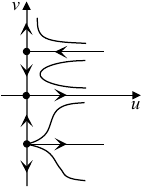}}~~~~~~
\subcaptionbox{%
     }{\includegraphics[height=1in]{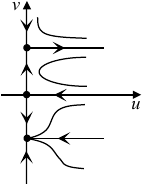}}~~~~~~
     \subcaptionbox{%
     }{\includegraphics[height=1in]{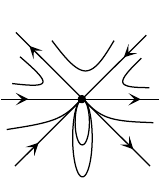}}
\caption{Phase portraits of (a) system~\eqref{equ:OIII+} and (b) system~\eqref{equ:OIII-}
along the $v$-axis, and (c) system~\eqref{GL:pfg} near $O$.}
\label{fig:OIII}
\end{figure}

Note that the existence and qualitative properties of equilibria on the $v$-axis
for system~\eqref{equ:OIII+} and \eqref{equ:OIII-} depend on the parties of $p,q,r$ and
also the signs of $b_q$ and $c_r$.
Then there are 32 situations.
However,
we only need to consider the following 12 situations:
\begin{center}
\begin{tabular}{lll}
{\bf(R1)} odd $p,q,r$, and $b_q,c_r>0$,
&
{\bf(R2)} odd $p,q,r$, $b_q>0$, $c_r<0$,
\\
{\bf(R3)} odd $p,q$, even $r$, $b_q, c_r>0$,
&
{\bf(R4)} odd $p,r$, even $q$,  $b_q,c_r>0$,
\\
{\bf(R5)} odd $p,r$, even $q$,  $b_q>0$, $c_r<0$,
&
{\bf(R6)} odd $p$, even $q,r$, $b_q,c_r>0$,
\\
{\bf(R7)} even $p$, odd $q,r$, $b_q,c_r>0$,
&
{\bf(R8)} even $p,r$, odd $q$, $b_q,c_r>0$,
\\
{\bf(R9)} even $p,q$, odd $r$, $b_q,c_r>0$,
&
{\bf(R10)} even $p,q$, odd $r$, $b_q<0$, $c_r>0$,
\\
{\bf(R11)} even $p,q,r$, $b_q, c_r>0$,
&
{\bf(R12)} even $p,q,r$, $b_q<0$, $c_r>0$,
\end{tabular}
\end{center}
since the other situations can be converted into the above situations
by the time-reversing and the transformation $x\to -x$,
or $y\to -y$, or $(x,y)\to (-x,-y)$.
In what follows,
we only give the proof in {\bf(R10)} since the proof in other situations are similar.
Hence,
the phase portraits of systems~\eqref{equ:OIII+} and \eqref{equ:OIII-}
along the $v$-axis are given by Fig.~\ref{fig:OIII} (a) and (b), respectively.
After blowing down,
we obtain the phase portrait Fig.~\ref{fig:OIII} (c),
or equivalently Fig.~\ref{fig:O-phase} (e),
of system~\eqref{GL:pfg} near the equilibrium $O$.
Consequently,
the proof of Theorem~\ref{thm:O} is completed.
\qquad$\Box$

\section{Phase portraits at the infinity}
\setcounter{equation}{0}
\setcounter{lm}{0}
\setcounter{thm}{0}
\setcounter{rmk}{0}
\setcounter{df}{0}
\setcounter{cor}{0}
\setcounter{pro}{0}

Note that system~\eqref{GL:pfg} with $\ell=1$ is equivalent to system~\eqref{equ:DH},
whose phase portraits at the infinity have been studied in \cite{DH99}.
So we further investigate system~\eqref{GL:pfg} with $\ell\ge 2$ near infinity.
There are 7 cases
{\small
\begin{center}
\begin{tabular}{llllll}
{\bf(C1)} $n=m=\ell$,
&{\bf(C2)} $n=m>\ell$,
&{\bf(C3)} $n=\ell>m$,
&{\bf(C4)} $n>\max\{\ell,m\}$,
\\
&{\bf(C5)} $m=\ell>n$,
&{\bf(C6)} $m>\max\{\ell,n\}$,
&{\bf(C7)} $\ell>\max\{m,n\}$.
\end{tabular}
\end{center}
}
\noindent
Due to the time-rescaling $t\to -t/|a_\ell|$,
we can assume without loss of generality that $a_\ell=-1$.
For convenience, we introduce the following notations:
\begin{center}
\begin{tabular}{lll}
{\bf(S1)} odd $\ell,m,n$, and $b_m,c_n>0$,
&
{\bf(S2)} odd $\ell,m,n$, $b_m>0$, $c_n<0$,
\\
{\bf(S3)} odd $\ell,m$, even $n$, $b_m, c_n>0$,
&
{\bf(S4)} odd $\ell,n$, even $m$, $b_m,c_n>0$,
\\
{\bf(S5)} odd $\ell,n$, even $m$, $b_m>0$, $c_n<0$,
&
{\bf(S6)} odd $\ell$, even $m,n$, $b_m,c_n>0$,
\\
{\bf(S7)} even $\ell$, odd $m,n$, $b_m,c_n>0$,
&
{\bf(S8)} even $\ell,n$, odd $m$, $b_m,c_n>0$,
\\
{\bf(S9)} even $\ell,m$, odd $n$, $b_m,c_n>0$,
&
{\bf(S10)} even $\ell,m$, odd $n$, $b_m<0$, $c_n>0$,
\\
{\bf(S11)} even $\ell,m,n$, $b_m, c_n>0$,
&
{\bf(S12)} even $\ell,m,n$, $b_m<0$, $c_n>0$,
\\
{\bf(T1)} odd $\ell,n$, $c_n>0$,
&{\bf(T2)} odd $\ell,n$, $c_n<0$,
\\
{\bf(T3)} odd $\ell$, even $n$, $c_n>0$,
&{\bf(T4)} even $\ell$, odd $n$, $c_n>0$,
\\
{\bf(T5)} even $\ell,n$, $c_n>0$.
\end{tabular}
\end{center}

\begin{figure}[H]
\centering
\subcaptionbox{%
     }{\includegraphics[height=1in]{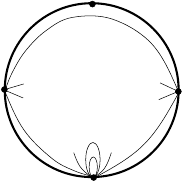}}~
\subcaptionbox{%
     }{\includegraphics[height=1in]{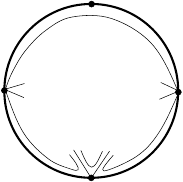}}~
\subcaptionbox{%
     }{\includegraphics[height=1in]{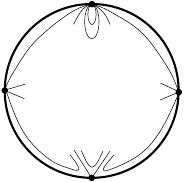}}~
\subcaptionbox{%
     }{\includegraphics[height=1in]{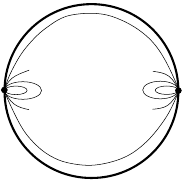}}~
\subcaptionbox{%
     }{\includegraphics[height=1in]{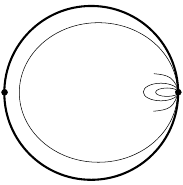}}
\\
\subcaptionbox{%
     }{\includegraphics[height=1in]{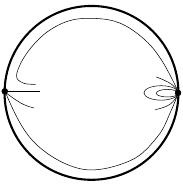}}~
\subcaptionbox{%
     }{\includegraphics[height=1in]{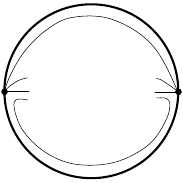}}~
\subcaptionbox{%
     }{\includegraphics[height=1in]{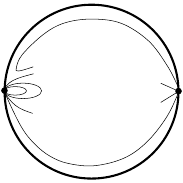}}~
\subcaptionbox{%
     }{\includegraphics[height=1in]{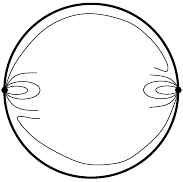}}~
\subcaptionbox{%
     }{\includegraphics[height=1in]{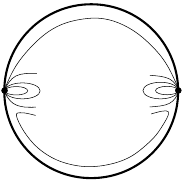}}
\\
\subcaptionbox{%
     }{\includegraphics[height=1in]{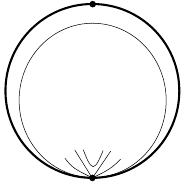}}~
\subcaptionbox{%
     }{\includegraphics[height=1in]{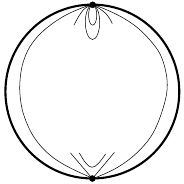}}~
\subcaptionbox{%
     }{\includegraphics[height=1in]{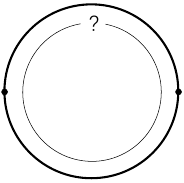}}~
\subcaptionbox{%
     }{\includegraphics[height=1in]{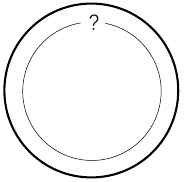}}
\caption{Topological phase portraits of system~\eqref{GL:pfg} at infinity.}
\label{fig:infty}
\end{figure}

\begin{thm}
System~\eqref{GL:pfg} with $a_\ell=-1$ and $\ell\ge 2$
has 26 different topological phase portraits near the equator of the Poincar\'{e} disc,
see Fig.~\ref{fig:infty}.
More concretely,
the relation between the phase portraits and parameters is listed in Table~\ref{tab:infty},
in which
$c_*:=\ell\big|\frac{b_m}{\ell+1}\big|^{\frac{\ell+1}{\ell}}$
and
$c^*:=\frac{(m-\ell)(n+1)}{(n-m)(\ell+1)}
\big|\frac{b_m(n-m)}{n-\ell}\big|^{\frac{n-\ell}{m-\ell}}$.
\label{thm:infty}
\end{thm}

\begin{center}
\small
\begin{longtable}{c|l|l}
\hline
\multirow{4}{*}{\bf(C1)}&
&{\bf(S1)} for Fig.\ref{fig:infty}(a),
\\
&&{\bf(S2)} with $c_n>-c_*$, $=-c_*$, $<-c_*$  for Fig.\ref{fig:infty}(b), (c), (x), resp.,
\\
&&{\bf(S11)} for Fig.\ref{fig:infty}(d),
\\
&&{\bf(S12)} with $c_n>c_*$, $=c_*$, $<c_*$ for Fig.\ref{fig:infty}(d), (e), (f), resp.
\\
\hline
\multirow{2}{*}{\bf(C2)}&
&{\bf(S1)} for Fig.\ref{fig:infty}(a),
{\bf(S2)} for Fig.\ref{fig:infty}(b),
{\bf(S6)} for Fig.\ref{fig:infty}(g),
\\
&&{\bf(S7)} for Fig.\ref{fig:infty}(h),
{\bf(S11)} for Fig.\ref{fig:infty}(i),
{\bf(S12)} for Fig.\ref{fig:infty}(e)
\\
\hline
{\bf(C3)}&
&
{\bf(T1)} for Fig.\ref{fig:infty}(a),
{\bf(T2)} for Fig.\ref{fig:infty}(x),
{\bf(T5)} for Fig.\ref{fig:infty}(d)
\\
\hline
\multirow{8}{*}{\bf(C4)}
&\multirow{2}{*}{$\frac{\ell(n+1)}{(\ell+1)m}>1$}
&
{\bf(T1)} for Fig.\ref{fig:infty}(n),
{\bf(T2)} for Fig.\ref{fig:infty}(w),
{\bf(T3)} for Fig.\ref{fig:infty}(d),
\\
&&
{\bf(T4)} for Fig.\ref{fig:infty}(o),
{\bf(T5)} for Fig.\ref{fig:infty}(d)
\\
\cline{2-3}
&$m=2\ell$
&
{\bf(S4)} for Fig.\ref{fig:infty}(n),
\\
&&
{\bf(S5)} with $c_n\ge -c^*$, $<-c^*$ for Fig.\ref{fig:infty}(u), (w), resp.,
\\
&
$n=2\ell+1$
&
{\bf(S9)} for Fig.\ref{fig:infty}(o),
\\
&&
{\bf(S10)} with $c_n> c^*$, $\le c^*$ for Fig.\ref{fig:infty}(o), (v), resp.
\\
\cline{2-3}
&$\frac{\ell(n+1)}{(\ell+1)m}=1$
&
{\bf(S1)} for Fig.\ref{fig:infty}(n),
\\
&$n>m+1$
&
{\bf(S2)} with $c_n\ge -c^*$, $<-c^*$ for Fig.\ref{fig:infty}(c), (w), resp.,
\\
&
&
{\bf(S4)} for Fig.\ref{fig:infty}(n),
\\
&&
{\bf(S5)} with $c_n\ge -c^*$, $<-c^*$ for Fig.\ref{fig:infty}(u), (w), resp.,
\\
&&
{\bf(S9)} for Fig.\ref{fig:infty}(o),
{\bf(S11)} for Fig.\ref{fig:infty}(d),
{\bf(S12)} for Fig.\ref{fig:infty}(d)
\\
\cline{2-3}
&\multirow{4}{*}{$\frac{\ell(n+1)}{(\ell+1)m}<1$}
&
{\bf(S1)} for Fig.\ref{fig:infty}(n),
{\bf(S2)} for Fig.\ref{fig:infty}(c),
{\bf(S3)} for Fig.\ref{fig:infty}(d),
\\
&&
{\bf(S4)} for Fig.\ref{fig:infty}(n),
{\bf(S5)} for Fig.\ref{fig:infty}(u),
{\bf(S6)} for Fig.\ref{fig:infty}(d),
\\
&&
{\bf(S7)} for Fig.\ref{fig:infty}(p),
{\bf(S8)} for Fig.\ref{fig:infty}(d),
{\bf(S9)} for Fig.\ref{fig:infty}(o),
\\
&&
{\bf(S10)} for Fig.\ref{fig:infty}(v),
{\bf(S11)} for Fig.\ref{fig:infty}(d),
{\bf(S12)} for Fig.\ref{fig:infty}(d)
\\
\hline
\multirow{3}{*}{\bf(C5)}&
&
{\bf(S1)} for Fig.\ref{fig:infty}(a),
{\bf(S2)} for Fig.\ref{fig:infty}(b),
{\bf(S3)} for Fig.\ref{fig:infty}(h),
\\
&&
{\bf(S9)} for Fig.\ref{fig:infty}(d),
{\bf(S10)} for Fig.\ref{fig:infty}(f),
{\bf(S11)} for Fig.\ref{fig:infty}(d),
\\
&&
{\bf(S12)} for Fig.\ref{fig:infty}(j)
\\
\hline
\multirow{4}{*}{\bf(C6)}&
&
{\bf(S1)} for Fig.\ref{fig:infty}(a),
{\bf(S2)} for Fig.\ref{fig:infty}(b),
{\bf(S3)} for Fig.\ref{fig:infty}(h),
\\
&&
{\bf(S4)} for Fig.\ref{fig:infty}(k),
{\bf(S5)} for Fig.\ref{fig:infty}(l),
{\bf(S6)} for Fig.\ref{fig:infty}(g),
\\
&&
{\bf(S7)} for Fig.\ref{fig:infty}(h),
{\bf(S8)} for Fig.\ref{fig:infty}(h),
{\bf(S9)} for Fig.\ref{fig:infty}(i),
\\
&&
{\bf(S10)} for Fig.\ref{fig:infty}(m),
{\bf(S11)} for Fig.\ref{fig:infty}(i),
{\bf(S12)} for Fig.\ref{fig:infty}(e)
\\
\hline
\multirow{8}{*}{\bf(C7)}
& \multirow{2}{*}{$m\le n$}
&
{\bf(T1)} for Fig.\ref{fig:infty}(n),
{\bf(T2)} for Fig.\ref{fig:infty}(w),
{\bf(T3)} for Fig.\ref{fig:infty}(o),
\\
&&
{\bf(T4)} for Fig.\ref{fig:infty}(d),
{\bf(T5)} for Fig.\ref{fig:infty}(d)
\\
\cline{2-3}
& \multirow{2}{*}{$m=n+1$}
&
{\bf(S3)} for Fig.\ref{fig:infty}(p),
{\bf(S4)} for Fig.\ref{fig:infty}(n),
{\bf(S5)} for Fig.\ref{fig:infty}(q),
\\
&&
{\bf(S8)} for Fig.\ref{fig:infty}(r),
{\bf(S9)} for Fig.\ref{fig:infty}(d),
{\bf(S10)} for Fig.\ref{fig:infty}(s)
\\
\cline{2-3}
&
\multirow{4}{*}{$m>n+1$}
&
{\bf(S1)} for Fig.\ref{fig:infty}(n),
{\bf(S2)} for Fig.\ref{fig:infty}(c),
{\bf(S3)} for Fig.\ref{fig:infty}(p),
\\
&&
{\bf(S4)} for Fig.\ref{fig:infty}(n),
{\bf(S5)} for Fig.\ref{fig:infty}(q),
{\bf(S6)} for Fig.\ref{fig:infty}(p),
\\
&&
{\bf(S7)} for Fig.\ref{fig:infty}(r),
{\bf(S8)} for Fig.\ref{fig:infty}(r),
{\bf(S9)} for Fig.\ref{fig:infty}(d),
\\
&&
{\bf(S10)} for Fig.\ref{fig:infty}(t),
{\bf(S11)} for Fig.\ref{fig:infty}(d),
{\bf(S12)} for Fig.\ref{fig:infty}(s)
\\
\hline
\caption{Parameter conditions for different topological phase portraits.}
\label{tab:infty}
\end{longtable}
\end{center}

\vspace{-50pt}
\begin{figure}[H]
\centering
\subcaptionbox{%
     }{\includegraphics[height=1in]{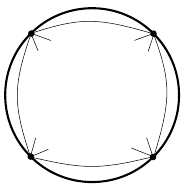}}~
\subcaptionbox{%
     }{\includegraphics[height=1in]{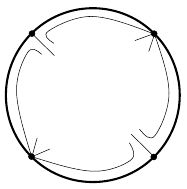}}~
\subcaptionbox{%
     }{\includegraphics[height=1in]{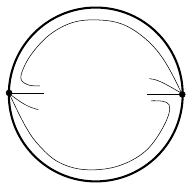}}~
\subcaptionbox{%
     }{\includegraphics[height=1in]{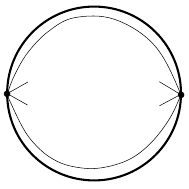}}~
\subcaptionbox{%
     }{\includegraphics[height=1in]{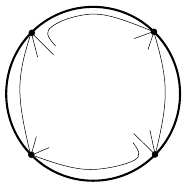}}
\\
\subcaptionbox{%
     }{\includegraphics[height=1in]{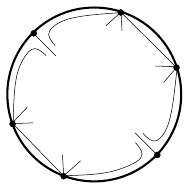}}~
\subcaptionbox{%
     }{\includegraphics[height=1in]{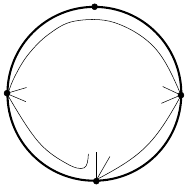}}~
\subcaptionbox{%
     }{\includegraphics[height=1in]{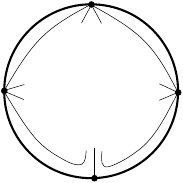}}~
\subcaptionbox{%
     }{\includegraphics[height=1in]{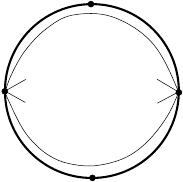}}~
\subcaptionbox{%
     }{\includegraphics[height=1in]{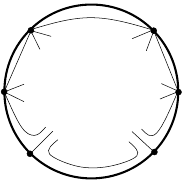}}
\end{figure}

\begin{rmk}
Any situation not lying in Table~\ref{tab:infty} can be converted to
a corresponding one in Table~\ref{tab:infty} by a time-reversing and
the transformation $x\to -x$, or $y\to -y$, or $(x,y)\to (-x,-y)$.
\label{rmk:otherS}
\end{rmk}

\begin{rmk}
System~\eqref{GL:pfg} with $a_\ell=-1$ and $\ell\ge 2$ is monodromic at infinity
if and only if either
\\
{\bf(W1)} $\ell(n+1)>(\ell+1)m$, $m\le n$, both $\ell$ and $n$ are odd and $c_n<0$, or
\\
{\bf(W2)} $\ell=m=n$ is odd, $b_m>0$ and $c_n<-c_*$, or
\\
{\bf(W3)} $\ell(n+1)=(\ell+1)m$, $n\ge m+1$,
both $\ell$ and $n$ are odd, and $c_n<-c^*$.
\label{rmk:MI}
\end{rmk}

{\bf Proof of Theorem~\ref{thm:infty}.}
Under the Poincar\'{e} transformation $x=1/v$ and $y=u/v$,
we obtain
\begin{equation}
\dot u=
-\frac{{\cal C}(v)}{v^{n-1}}
-\frac{u{\cal A}(u,v)}{v^{\ell-1}}
+\frac{u{\cal B}(v)}{v^{m-1}},~~~
\dot v=
-\frac{{\cal A}(u,v)}{v^{\ell-2}}
+\frac{{\cal B}(v)}{v^{m-2}},
\label{equ:8C1}
\end{equation}
where
$$
{\cal A}(u,v):=\sum_{i=0}^\ell a_{\ell-i} u^{\ell-i}v^i,~~~
{\cal B}(v):=\sum_{i=0}^m b_{m-i} v^{i},~~~
{\cal C}(v):=\sum_{i=0}^n c_{n-i} v^{i}.
$$
On the other hand,
under the Poincar\'{e} transformation $x=u/v$ and $y=1/v$,
we obtain
\begin{equation}
\dot u=
\frac{\widetilde{\cal A}(v)}{v^{\ell-1}}
-\frac{\widetilde{\cal B}(u,v)}{v^{m-1}}
+\frac{u\widetilde{\cal C}(u,v)}{v^{n-1}},~~~
\dot v=
\frac{\widetilde{\cal C}(u,v)}{v^{n-2}},
\label{equ:8C2}
\end{equation}
where
$$
\widetilde{\cal A}(v):=\sum_{i=0}^\ell a_{\ell-i} v^i,~~~
\widetilde{\cal B}(u,v):=\sum_{i=0}^m b_{m-i} u^{m-i} v^i,~~~
\widetilde{\cal C}(u,v):=\sum_{i=0}^n c_{n-i} u^{n-i} v^i.
$$
In order to cancel factors of $v$ in denominators,
we need to multiply systems~\eqref{equ:8C1} and \eqref{equ:8C2} by $v^{\max\{\ell,m,n\}-1}$
and therefore there are 7 cases {\bf(C1)}-{\bf(C7)},
given just before Theorem~\ref{thm:infty}.
We only give proofs in cases {\bf(C1)}, {\bf(C2)} and {\bf(C7)}
since cases {\bf(C3)} and {\bf(C5)} are similar to {\bf(C2)}
and cases {\bf(C4)} and {\bf(C6)} are similar to {\bf(C7)}.

\begin{figure}[h]
\centering
\subcaptionbox{%
     }{\includegraphics[height=1in]{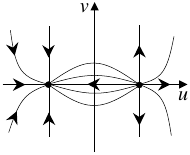}}~~
\subcaptionbox{%
     }{\includegraphics[height=1in]{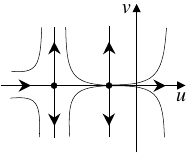}}~~
\subcaptionbox{%
     }{\includegraphics[height=1in]{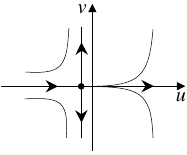}}~~
\subcaptionbox{%
     }{\includegraphics[height=1in]{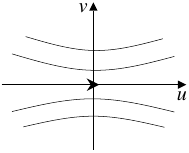}}
\\
\subcaptionbox{%
     }{\includegraphics[height=1in]{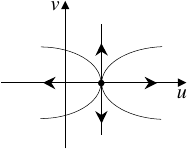}}~~
\subcaptionbox{%
     }{\includegraphics[height=1in]{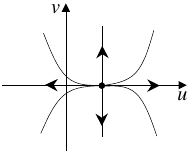}}~~
\subcaptionbox{%
     }{\includegraphics[height=1in]{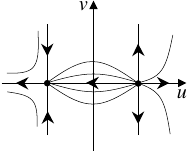}}~~
\subcaptionbox{%
     }{\includegraphics[height=1in]{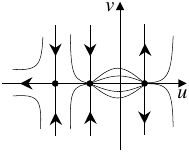}}
\caption{Phase portraits in {\bf(C1)}.
}
\label{fig:C1}
\end{figure}

In case {\bf(C1)},
under the time rescaling $dt\to v^{n-1}dt$,
system~\eqref{equ:8C1} becomes
\begin{equation}
\dot u=E(u)+O(v),~~~
\dot v=v(b_m+ u^\ell+O(v)),
\label{equ:8C1t}
\end{equation}
where $E(u):=-c_n+b_m u+ u^{\ell+1}$,
and system~\eqref{equ:8C2} becomes
\begin{equation}
\dot u=-1+o(1),~~~\dot v=O(v)
\label{equ:8C2t}
\end{equation}
near the origin $(0,0)$.
Clearly,
the origin is a regular point of system~\eqref{equ:8C2t}.
If $(u_*,0)$ is an equilibrium of system~\eqref{equ:8C1t} on the $u$-axis,
then the Jacobian matrix at this equilibrium is
\begin{equation}
\left(
\begin{array}{ccc}
E'(u_*)  &  \star
\\
0        &  c_n/u_*
\end{array}
\right),
\label{equ:JIC1}
\end{equation}
where we used the fact that $E(u_*)=0$.
In case {\bf(C1)},
we only need to consider situations {\bf(S1)}, {\bf(S2)}, {\bf(S11)} and {\bf(S12)},
given just before Theorem~\ref{thm:infty}.
For {\bf(S1)},
the derivative $E'$ has exactly one (negative) real zero $u=(\frac{-b_m}{\ell+1})^{1/\ell}$
and therefore $E$ has one positive zero $u_1$ and one negative zero $u_2$ since $E(0)=-c_n<0$.
Then we see from \eqref{equ:JIC1} that
the equilibrium $(u_1,0)$ (and $(u_2,0)$) is an unstable (and stable) node,
and the phase portrait of system~\eqref{equ:8C1t}
along the $u$-axis is given by Fig.~\ref{fig:C1}(a).
So we obtain the topologival phase portrait Fig.~\ref{fig:infty}(a)
of system~\eqref{GL:pfg} near the equator of the Poincar\'{e} disc.

For {\bf(S2)},
the derivative $E'$ has exactly one (negative) real zero $u=(\frac{-b_m}{\ell+1})^{1/\ell}$
and then $E$ has
two negative real zeros if $c_n>-c_*$,
one (double) real zero $u=(\frac{-b_m}{\ell+1})^{1/\ell}$ if $c_n=-c_*$,
and no real zeros if $c_n<-c_*$.
Qualitative property of each equilibrium can be obtained from \eqref{equ:JIC1} directly.
Then phase portraits of system~\eqref{equ:8C1t} along the $u$-axis are given by
Fig.~\ref{fig:C1}(b)-(d) for $c_n>-c_*$, $=-c_*$ and $<-c_*$, respectively.
Consequently,
we obtain topological phase portraits Fig.~\ref{fig:infty}(b), (c), and (x)
of system~\eqref{GL:pfg} near the equator of the Poincar\'{e} disc
for $c_n>-c_*$, $=-c_*$, and $<-c_n^*$, respectively.

For {\bf(S11)},
the derivative $E'$ has no real zeros and therefore
$E$ has exactly one (positive) real zero.
It follows from \eqref{equ:JIC1} that
the only equilibrium of system~\eqref{equ:8C1t} on the $u$-axis is an unstable node,
as illustrated by Fig.~\ref{fig:C1}(e).
from which we obtain the topological phase portrait Fig.~\ref{fig:infty}(d)
of system~\eqref{GL:pfg} near the equator of the Poincar\'{e} disc.

For {\bf(S12)},
the derivative $E'$ has two real zeros
$u_\pm:=\pm (\frac{-b_m}{\ell+1})^{1/\ell}$ and then
$E$ has exactly one real zero ($>u_+$) if $c_n>c_*=E(u_-)$,
a double zero $u_-$ and one positive zero if $c_n=c_*$,
two negative zeros and one positive zero if $c_n<c_*$.
Qualitative property of each equilibrium can be obtained from \eqref{equ:JIC1} directly.
Then phase portraits of system~\eqref{equ:8C1t} along the $u$-axis are given by Fig.~\ref{fig:C1}(f)-(h),
for $c_n>c_*$, $=c_*$ and $<c_*$, respectively.
So we obtain the corresponding topological phase portraits Fig.~\ref{fig:infty}(d)-(f)
of system~\eqref{GL:pfg} near the equator of the Poincar\'{e} disc.

\begin{figure}[H]
\centering
\includegraphics[height=1in]{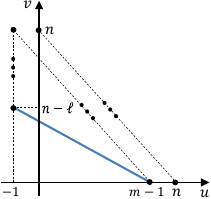}
\caption{Newton polygon of system~\eqref{equ:8C2C2t}.}
\label{fig:NPC2}
\end{figure}

In case {\bf(C2)},
under the time rescaling $dt\to v^{n-1}dt$,
system~\eqref{equ:8C1} becomes
\begin{eqnarray}
\dot u=-c_n+b_mu+O(v),~~~
\dot v=v(b_m+O(v)),
\label{8C1C2t}
\end{eqnarray}
and system~\eqref{equ:8C2} becomes
\begin{equation}
\dot u=v^{n-\ell}\widetilde{\cal A}(v)-\widetilde{\cal B}(u,v)+u\widetilde{\cal C}(u,v),~~~
\dot v=v\widetilde{\cal C}(u,v),
\label{equ:8C2C2t}
\end{equation}
Clearly,
$(c_n/b_m,0)$ is the only equilibrium of system~\eqref{8C1C2t} on the $u$-axis,
which is a stable (or unstable) node if $b_m<0$ (or $>0$).
The equilibrium $(0,0)$ of system~\eqref{equ:8C2C2t} is non-elementary and
we need to desingularize this equilibrium to obtain its qualitative property.
Note that each support point of system~\eqref{equ:8C2C2t} lies on either the line $u=-1$,
or the line $u+v=m-1$, or the line $u+v=n$, as illustrated by Fig.~\ref{fig:NPC2}.
Then the Newton polygon has only one edge,
lying on the line $\alpha_1 u+\beta_1 v=\delta_1:=(n-\ell)(m-1)$,
where $\alpha_1:=n-\ell$ and $\beta_1:=m$.
We rewrite system~\eqref{equ:8C2C2t}
in quasi-homogeneous components of type $(\alpha_1,\beta_1)$ as
\begin{equation*}
\dot u=\sum_{k=\delta_1}^{\delta_*}P_k^{(1)}(u,v)+\cdots,~~~
\dot v=\sum_{k=\delta_1}^{\delta_*}Q_k^{(1)}(u,v)+\cdots
\end{equation*}
such that
$P_k^{(1)}(r^{\alpha_1}u,r^{\alpha_1}v)=r^{\alpha_1+k}P_k^{(1)}(u,v)$ and
$Q_k^{(1)}(r^{\alpha_1}u,r^{\alpha_1}v)=r^{\beta_1+k}Q_k^{(1)}(u,v)$,
where
$P_{\delta_1}^{(1)}(u,v):=-v^{n-\ell}-b_mu^m$,
$Q_{\delta_1}^{(1)}(u,v)=0$, $\delta_*:=\alpha_1 n=\delta_1+(n-\ell)$,
and dots represent those terms of quasi-homogeneous degree bigger than $\delta_*$.
Note that
all support points of $\dot v$ in system~\eqref{equ:8C2C2t} lie on the line $u+v=n$,
on which $(n,0)$ is the only support point $(i,j)$ such that
the sum $\alpha_1 i+\beta_1 j$ takes its minimum $\delta_*$.
Then
$$
Q_{\delta_1+1}^{(1)}(u,v)=\cdots=Q_{\delta_*-1}^{(1)}(u,v)=0~~~\mbox{and}~~~
Q_{\delta_*}^{(1)}(u,v)=c_nu^nv.
$$
So blowing up the equilibrium $(0,0)$ of system~\eqref{equ:8C2C2t}
in the positive $u$-direction
by the transformation $u=u_1^{\alpha_1}$ and $v=u_1^{\beta_1}v_1$,
we obtain
\begin{equation}
\dot u_1={\cal U}_1(u_1,v_1),~~~
\dot v_1={\cal V}_1(u_1,v_1),
\label{equ:8C2C2tb}
\end{equation}
where a time rescaling is performed and
{\small
\begin{align*}
{\cal U}_1(u_1,v_1)
&:=u_1\{W_1(u_1,v_1)+u_1^{\delta_*-\delta_1}P_{\delta_*}^{(1)}(1,v_1)
+O(u_1^{\delta_*-\delta_1+1})\},
\\
{\cal V}_1(u_1,v_1)
&:=-\beta_1 v_1 W_1(u_1,v_1)
+u_1^{\delta_*-\delta_1}\{\alpha_1Q_{\delta_*}^{(1)}(1,v_1)
-\beta_1 v_1P_{\delta_*}^{(1)}(1,v_1)\}
+O(u_1^{\delta_*-\delta_1+1}),
\\
W_1(u_1,v_1)
&:=P_{\delta_1}^{(1)}(1,v_1)+u_1P_{\delta_1+1}^{(1)}(1,v_1)+\cdots
+u_1^{\delta_*-\delta_1-1}P_{\delta_*-1}^{(1)}(1,v_1).
\end{align*}
}On the $v_1$-axis,
system~\eqref{equ:8C2C2tb} has the equilibrium $(0,0)$ and
at most two equilibria determined by the equation
$$
W_1(0,v_1)=P_{\delta_1}^{(1)}(1,v_1)=-v_1^{n-\ell}-b_m=0.
$$
Jacobian matrices of system~\eqref{equ:8C2C2tb} at the equilibrium $(0,0)$ and the equilibrium $(0,v_1^*)$ with $v_1^*\ne 0$ (if exists) are given by
\begin{equation*}
\left(
\begin{array}{cccc}
-b_m  & 0
\\
\star & \beta_1 b_m
\end{array}
\right)~~~\mbox{and}~~~
\left(
\begin{array}{cccc}
0     & 0
\\
\star & \beta_1(n-\ell)(v_1^*)^{n-\ell}
\end{array}
\right),
\end{equation*}
respectively.
Clearly, the equilibrium $(0,0)$ is a hyperbolic saddle and
the equilibrium $(0,v_1^*)$ (if exists) is semi-hyperbolic.
In order to determine the qualitative property,
by Theorem~7.1 in \cite[p.114]{ZZF},
we solve from the equation ${\cal V}_1(u_1,v_1)=0$
near the semi-hyperbolic equilibrium $(0,v_1^*)$
that $v_1=\Lambda_1(u_1)=v_1^*+O(u_1)$ and
$$
W_1(u_1,\Lambda_1(u_1))=\frac{\alpha_1}{\beta_1} c_n u_1^{\delta_*-\delta_1}
-P_{\delta_*}^{(1)}(1,v_1)u_1^{\delta_*-\delta_1}
+O(u_1^{\delta_*-\delta_1+1}).
$$
Substituting it in ${\cal U}_1(u_1,v_1)$ of \eqref{equ:8C2C2tb}, we obtain that
$$
{\cal U}_1(u_1,\Lambda_1(u_1))
=\frac{\alpha_1}{\beta_1}c_n u_1^{\delta_*-\delta_1+1}+O(u_1^{\delta_*-\delta_1+2}).
$$
Then the qualitative property of the semi-hyperbolic equilibrium $(0,v_1^*)$ (if exists)
of system~\eqref{equ:8C2C2tb} is determined by
the parity of $\delta_*-\delta_1+1$
and the signs of $\beta_1(n-\ell)(v_1^*)^{n-\ell}$ and $\alpha_1c_n/\beta_1$.

On the other hand,
blowing up the equilibrium $(0,0)$ of system~\eqref{equ:8C2C2t}
in the negative $u$-direction
by the transformation $u=-u_1^{\alpha_1}$ and $v=u_1^{\beta_1}v_1$, we obtain
\begin{equation}
\dot u_1=\widetilde{\cal U}_1(u_1,v_1),~~~
\dot v_1=\widetilde{\cal V}_1(u_1,v_1),
\label{equ:8C2C2tb-}
\end{equation}
where a time rescaling is performed and
{\small
\begin{align*}
\widetilde{\cal U}_1(u_1,v_1)
&:=-u_1\{\widetilde{W}_1(u_1,v_1)+u_1^{\delta_*-\delta_1}P_{\delta_*}^{(1)}(-1,v_1)
+O(u_1^{\delta_*-\delta_1+1})\},
\\
\widetilde{\cal V}_1(u_1,v_1)
&:=\beta_1 v_1 \widetilde{W}_1(u_1,v_1)
+u_1^{\delta_*-\delta_1}\{\alpha_1Q_{\delta_*}^{(1)}(-1,v_1)
+\beta_1 v_1 P_{\delta_*}^{(1)}(-1,v_1)\}
+O(u_1^{\delta_*-\delta_1+1}),
\\
\widetilde{W}_1(u_1,v_1)
&:=P_{\delta_1}^{(1)}(-1,v_1)+u_1P_{\delta_1+1}^{(1)}(-1,v_1)+\cdots
+u_1^{\delta_*-\delta_1}P_{\delta_*-1}^{(1)}(-1,v_1).
\end{align*}
}System~\eqref{equ:8C2C2tb-} has the equilibrium $(0,0)$ and at most two equilibria determined by the equation
$$
\widetilde{W}_1(0,v_1)=P_{\delta_1}^{(1)}(-1,v_1)=-v_1^{n-\ell}-(-1)^mb_m=0.
$$
Jacobian matrices of system~\eqref{equ:8C2C2tb-} at the equilibrium $(0,0)$ and the equilibrium $(0,\tilde{v}_1^*)$ with $\tilde{v}_1^*\ne 0$ (if exists) are given by
\begin{equation*}
\left(
\begin{array}{cccc}
(-1)^mb_m & 0
\\
\star     & -\beta_1(-1)^mb_m
\end{array}
\right)~~~\mbox{and}~~~
\left(
\begin{array}{cccc}
0     & 0
\\
\star & -\beta_1(n-\ell)(\tilde{v}_1^*)^{n-\ell}
\end{array}
\right),
\end{equation*}
respectively.
Clearly, the equilibrium $(0,0)$ is a hyperbolic saddle and
the equilibrium $(0,\tilde{v}_1^*)$ (if exists) is semi-hyperbolic.
Similar to the equilibrium $(0,\tilde{v}_1^*)$ of system~\eqref{equ:8C2C2tb},
qualitative property of the semi-hyperbolic equilibrium $(0,\tilde{v}_1^*)$ is determined by the parity of $\delta_*-\delta_1+1$ and
the signs of $-\beta_1(n-\ell)(\tilde{v}_1^*)^{n-\ell}$ and $(-1)^n\alpha_1c_n/\beta_1$.

\begin{figure}[H]
\centering
\subcaptionbox{%
     }{\includegraphics[height=1in]{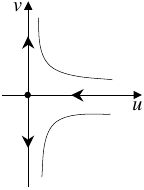}}~~
\subcaptionbox{%
     }{\includegraphics[height=1in]{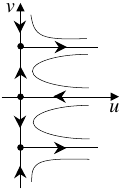}}~~
\subcaptionbox{%
     }{\includegraphics[height=1in]{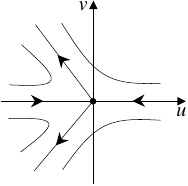}}
\caption{Phase portraits in {\bf(C2)}.}
\label{fig:C2}
\end{figure}

Note that in case {\bf(C2)} we only need to consider situations
{\bf(S1)}, {\bf(S2)}, {\bf(S6)}, {\bf(S7)}, {\bf(S11)} and {\bf(S12)}.
We only give the proof in {\bf(S2)} and the other situations can be discussed similarly.
For {\bf(S2)},
system~\eqref{equ:8C2C2tb} has the only equilibrium $(0,0)$ on the $v_1$-axis
and the phase portrait along the $v_1$-axis is given by Fig.~\ref{fig:C2}(a).
System~\eqref{equ:8C2C2tb-} has three equilibria on the $v_1$-axis
and the phase portrait along the $v_1$-axis is given by Fig.~\ref{fig:C2}(b).
After blowing down,
we obtain the phase portrait Fig.~\ref{fig:C2}(c) of system~\eqref{equ:8C2C2t}
at the origin.
Note that the only equilibrium of system~\eqref{8C1C2t} on the $u$-axis is an unstable node.
Then the topological phase portrait of system~\eqref{GL:pfg}
near the equator of the Poincar\'{e} disc is given by Fig.~\ref{fig:infty}(b).

\begin{figure}[H]
\centering
\subcaptionbox{%
     }{\includegraphics[height=1in]{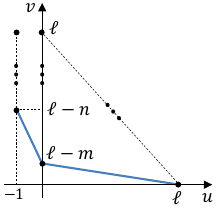}}~~~~~~
\subcaptionbox{%
     }{\includegraphics[height=1in]{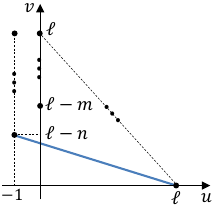}}
\caption{Newton polygon of system~\eqref{equ:8C1C7t}.
(a) for {\bf (C7a)} and (b) for {\bf (C7b)}.}
\label{fig:NP-C7}
\end{figure}

In case {\bf(C7)},
under the time rescaling $dt\to v^{\ell-1}dt$,
system~\eqref{equ:8C1} becomes
\begin{eqnarray}
\dot u=-v^{\ell-n}{\cal C}(v)-u{\cal A}(u,v)+uv^{\ell-m}{\cal B}(v),~~~
\dot v=-v{\cal A}(u,v)+v^{\ell-m+1}{\cal B}(v),
\label{equ:8C1C7t}
\end{eqnarray}
and system~\eqref{equ:8C2} becomes
\begin{equation}
\dot u=-1+O(v),~~~
\dot v=O(v)
\label{equ:8C2C7t}
\end{equation}
near the origin.
Clearly,
the origin is a regular point of system~\eqref{equ:8C2C7t}
and system~\eqref{equ:8C1C7t} has the only equilibrium $(0,0)$ on the $u$-axis,
which is non-elementary.
Each support point of system~\eqref{equ:8C1C7t} lies on either the line $u=-1$,
or the line $u=0$, or the line $u+v=\ell$, as illustrated by Fig.~\ref{fig:NP-C7}.
There are two circumstances
{\bf (C7a)} $m\le n$, and
{\bf (C7b)} $m>n$.

\begin{figure}[H]
\centering
\subcaptionbox{%
     }{\includegraphics[height=1in]{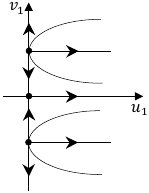}}~~
\subcaptionbox{%
     }{\includegraphics[height=1in]{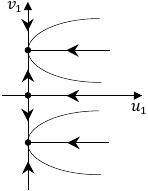}}~~
\subcaptionbox{%
     }{\includegraphics[height=1in]{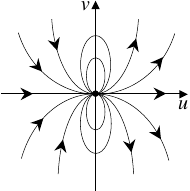}}~~
\caption{Phase portraits of systems~\eqref{equ:8C4ap},
\eqref{equ:8C4an} and \eqref{equ:8C1C7t}.}
\label{fig:C7aT1}
\end{figure}

In {\bf(C7a)},
the Newton polygon of system~\eqref{equ:8C1C7t} has only one edge,
linking vertex $(-1,\ell-n)$ with vertex $(\ell,0)$,
and there are no support points on this edge except for the two vertices,
see Fig.~\ref{fig:NP-C7}(a).
The only edge lies on the line $\alpha_2 u+\beta_2 v=\delta_2:=\ell(\ell-n)$,
where $\alpha_2:=\ell-n$ and $\beta_2:=\ell+1$.
We rewrite system~\eqref{equ:8C1C7t}
in quasi-homogeneous components of type $(\alpha_2,\beta_2)$ as
\begin{equation*}
\dot u=P_{\delta_2}^{(2)}(u,v)+\cdots,~~~
\dot v=Q_{\delta_2}^{(2)}(u,v)+\cdots,
\end{equation*}
where
$P_{\delta_2}^{(2)}(u,v):=-c_nv^{\ell-n}+u^{\ell+1}$,
$Q_{\delta_2}^{(2)}(u,v)=u^\ell v$,
and dots represent those terms of quasi-homogeneous degree bigger than $\delta_2$.
Blowing up the non-elementary equilibrium $(0,0)$ of system~\eqref{equ:8C1C7t}
in the positive $u$-direction
by the transformation $u=u_1^{\alpha_2}$ and $v=u_1^{\beta_2}v_1$, we obtain
\begin{eqnarray}
\dot u_1=u_1\{P_{\delta_2}^{(2)}(1,v_1)+O(u_1)\},~~~
\dot v_1=G_0(v_1)+O(u_1),
\label{equ:8C4ap}
\end{eqnarray}
where a time rescaling is performed and
$$
G_0(v_1):=\alpha_2 Q_{\delta_2}^{(2)}(1,v_1)-\beta_2 v_1 P_{\delta_2}^{(2)}(1,v_1)
=v_1\{\beta_2 c_n v_1^{\ell-n}-(\beta_2-\alpha_2)\}.
$$
On the other hand,
blowing up the equilibrium $(0,0)$ of system~\eqref{equ:8C1C7t}
in the negative $u$-direction
by the transformation $u=-u_1^{\alpha_2}$ and $v=u_1^{\beta_2}v_1$, we obtain
\begin{eqnarray}
\dot u_1=-u_1\{P_{\delta_2}^{(2)}(-1,v_1)+O(u_1)\},~~~
\dot v_1=\widetilde{G}_0(v_1)+O(u_1),
\label{equ:8C4an}
\end{eqnarray}
where a time rescaling is performed and
$$
\widetilde{G}_0(v_1)
:=\alpha_2Q_{\delta_2}^{(2)}(-1,v_1)+\beta_2v_1P_{\delta_2}^{(2)}(-1,v_1)
=-v_1\{\beta_2 c_n v_1^{\ell-n}+(-1)^{\ell}(\beta_2-\alpha_2)\}.
$$
Jacobian matrices of system~\eqref{equ:8C4ap} at the equilibrium $(0,0)$ and
the equilibrium $(0,v_1^*)$ with $v_1^*\ne 0$ (if exists) are given by
\begin{equation*}
\left(
\begin{array}{cccc}
1     & 0
\\
\star & -(n+1)
\end{array}
\right)~~~\mbox{and}~~~
\left(
\begin{array}{cccc}
\frac{\alpha_2}{\beta_2} & 0
\\
\star                    & G'_0(v_1^*)
\end{array}
\right),
\end{equation*}
respectively.
Jacobian matrices of system~\eqref{equ:8C4an} at the equilibrium $(0,0)$ and
the equilibrium $(0,\tilde{v}_1^*)$ with $\tilde{v}_1^*\ne 0$ (if exists) are given by
\begin{equation*}
\left(
\begin{array}{cccc}
(-1)^\ell & 0
\\
\star     & (-1)^{\ell+1}(n+1)
\end{array}
\right)~~~\mbox{and}~~~
\left(
\begin{array}{cccc}
(-1)^\ell\frac{\alpha_2}{\beta_2} & 0
\\
\star                             & \widetilde{G}'_0(\tilde{v}_1^*)
\end{array}
\right),
\end{equation*}
respectively.
Since the above analyses are independent of $m$ and $b_m$,
we only need to consider situations
{\bf(T1)}-{\bf(T5)},
given just before Theorem~\ref{thm:infty},
and we only give the proof for {\bf(T1)}
since proofs for other situations are similar.
For {\bf(T1)},
both systems~\eqref{equ:8C4ap} and \eqref{equ:8C4an} have three equilibria
on the $v_1$-axis, and their phase portraits along the $v_1$-axis
are given by Fig.~\ref{fig:C7aT1}(a) and (b), respectively.
After blowing down,
we obtain the phase portrait Fig.~\ref{fig:C7aT1}(c) of system~\eqref{equ:8C1C7t}
at the origin and then the phase portrait Fig~\ref{fig:infty}(n) of system~\eqref{GL:pfg}
near the equator of the Poincar\'{e} disc.

In {\bf (C7b)},
the Newton polygon has exactly two edges,
linking vertices $(-1,\ell-n)$, $(0,\ell-m)$ and $(\ell,0)$,
see Fig.~\ref{fig:NP-C7}(b).
The first edge lies on the line $\alpha_3 u+\beta_3 v=\delta_3:=\ell-m$,
where $\alpha_3:=m-n$ and $\beta_3:=1$.
We rewrite system~\eqref{equ:8C1C7t}
in quasi-homogeneous components of type $(\alpha_3,\beta_3)$ as
\begin{eqnarray}
\dot u=P_{\delta_3}^{(3)}(u,v)+\cdots,~~~
\dot v=Q_{\delta_3}^{(3)}(u,v)+\cdots,
\end{eqnarray}
where
$P_{\delta_3}^{(3)}(u,v):=-c_nv^{\ell-n}+b_muv^{\ell-m}$,
$Q_{\delta_3}^{(3)}(u,v):=b_mv^{\ell-m+1}$,
and dots represent those terms of quasi-homogeneous degree bigger than $\delta_3$.
Blowing up the non-elementary equilibrium $(0,0)$ of system~\eqref{equ:8C1C7t}
in the positive $u$-direction
by the transformation $u=u_1^{\alpha_3}$ and $v=u_1^{\beta_3}v_1$,
we obtain
\begin{equation}
\left\{
\begin{array}{llll}
\dot u_1=u_1\frac{1}{(u_1)^{\alpha_3+\delta_3}}\dot u
&=u_1\{P_{\delta_3}^{(3)}(1,v_1)+O(u_1)\},
\\
\dot v_1=\frac{\alpha_3}{(u_1)^{\beta_3+\delta_3}}\dot v -\frac{\beta_3v_1}{(u_1)^{\alpha_3+\delta_3}}\dot u
&=G_1(v_1)+O(u_1),
\end{array}
\right.
\label{equ:C7bp}
\end{equation}
where a time rescaling is performed and
$G_1(v_1):=v_1^{\ell-m+1}\{\beta_3c_nv_1^{m-n}+(\alpha_3-\beta_3)b_m\}$.
System~\eqref{equ:C7bp} has the equilibrium $(0,0)$ and
at most two other equilibria on the $v_1$-axis.
The equilibrium $(0,0)$ is non-elementary and
the Jacobian matrix at the equilibrium $(0,v_1^*)$ with $v_1^*\ne 0$ (if exists) is given by
\begin{equation}
\left(
\begin{array}{cccc}
\frac{\alpha_3}{\beta_3}b_m(v_1^*)^{\ell-m} & 0
\\
\star & G'_1(v_1^*)
\end{array}
\right).
\label{equ:JC7b21}
\end{equation}
In order to investigate the property of
the non-elementary equilibrium $(0,0)$ of system~\eqref{equ:C7bp},
we need the following fact.

\begin{figure}[H]
\centering
\includegraphics[height=1in]{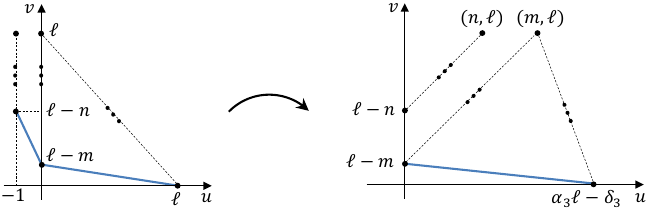}
\caption{Newton polygon of system~\eqref{equ:8C1C7t}.
(a) for {\bf (C7a)} and (b) for {\bf (C7b)}.}
\label{fig:C7NPBNP}
\end{figure}

\noindent
{\bf Fact~1.}
{\it The Newton polygon of system~\eqref{equ:C7bp} has only one edge,
linking vertex $(0,\alpha_4)$ with vertex $(\beta_4,0)$
and lying on the line $\alpha_4 u+\beta_4 v=\delta_4:=\alpha_4\beta_4$,
where $\alpha_4:=\ell-m$ and $\beta_4:=\alpha_3\ell-\delta_3$.
Moreover,
system~\eqref{equ:C7bp} can be rewritten
in quasi-homogeneous components of type $(\alpha_4,\beta_4)$ as
\begin{equation}
\dot u_1=\sum_{k=\delta_4}^{\delta'}P_k^{(4)}(u_1,v_1)+\cdots,~~~
\dot v_1=\sum_{k=\delta_4}^{\delta'}Q_k^{(4)}(u_1,v_1)+\cdots,
\end{equation}
where $\delta':=\beta_4(\ell-n)$,
dots represent those terms of quasi-homogeneous degree bigger than $\delta'$,
\begin{eqnarray*}
P_k^{(4)}(u_1,v_1)=u_1R_k(u_1,v_1),~~~
Q_k^{(4)}(u_1,v_1)=(\alpha_3-\beta_3)v_1R_k(u_1,v_1),
\end{eqnarray*}
for all $k=\delta_4,...,\delta'-1$,
\begin{eqnarray*}
&&P_{\delta'}^{(4)}(u_1,v_1)=u_1\{-c_nv_1^{\ell-n}+R_{\delta'}(u_1,v_1)\},
\\
&&Q_{\delta'}^{(4)}(u_1,v_1)=v_1\{\beta_3c_nv_1^{\ell-n}
+(\alpha_3-\beta_3)R_{\delta'}(u_1,v_1)\},
\end{eqnarray*}
and
$
R_k(u_1,v_1):=b_{\ell-j'_k}u_1^{i'_k}v_1^{j'_k}-a_{\ell-j''_k}u_1^{i''_k}v_1^{j''_k}
$
for all $k=\delta_4,...,\delta'$,
the point $(i'_k,j'_k)$ $($and $(i''_k,j''_k)$$)$ is the interception point of
the line $v=(u+\delta_3)/\beta_3$
$($and the line $v=(u-\alpha_3\ell+\delta_3)/(\beta_3-\alpha_3)$$)$ and
the line $\alpha_4 u+\beta_4 v=k$,
moreover,
$b_{\ell-j'_k}:=0$ $($and $a_{\ell-j''_k}:=0$$)$
if $j'_k$ $($and $j''_k$$)$ is not an integer.
Especially, $R_{\delta_4}^{(4)}(u_1,v_1)=b_mv_1^{\ell-m}+u_1^{\alpha_3\ell-\delta_3}$.
}

Actually,
we see from \eqref{equ:C7bp} that
each support point $(i,j)$ of system~\eqref{equ:8C1C7t} becomes a support point
$(\alpha_4i+\beta_4j-\delta_3,j)$ of system~\eqref{equ:C7bp}.
So all support points of system~\eqref{equ:8C1C7t} lying on the line $u=-1$,
the line $u=0$ and
the line $u+v=\ell$
become support points of system~\eqref{equ:C7bp}
lying on the line $v=(u+\alpha_3+\delta_3)/\beta_3$,
the line $v=(u+\delta_3)/\beta_3$, and
the line $v=(u-\alpha_3\ell+\delta_3)/(\beta_3-\alpha_3)$,
respectively.
Moreover,
each support point $(-1,j)$ of system~\eqref{equ:8C1C7t} has coefficient $(-c_{\ell-j},0)$
and becomes a support point $(\beta_3j-\alpha_3-\delta_3,j)$ of system~\eqref{equ:C7bp}
with coefficient $(-c_{\ell-j},\beta_3c_{\ell-j})$,
each support point $(0,j)$ of system~\eqref{equ:8C1C7t} has
coefficient $(b_{\ell-j},b_{\ell-j})$
and becomes a support point $(\beta_3j-\delta_3,j)$ of system~\eqref{equ:C7bp}
with coefficient $(b_{\ell-j},(\alpha_3-\beta_3)b_{\ell-j})$,
and
each support point $(\ell-j,j)$ of system~\eqref{equ:8C1C7t}
has coefficient $(-a_{\ell-j},-a_{\ell-j})$ and
becomes a support point $((\beta_3-\alpha_3)j+\alpha_3\ell-\delta_3,j)$
of system~\eqref{equ:C7bp} with coefficient
$(-a_{\ell-j},-(\alpha_3-\beta_3)a_{\ell-j})$.
So we obtain those quasi-homogeneous components of type $(\alpha_4,\beta_4)$
given in {\bf Fact~1}, which completes the proof.

By {\bf Fact~1},
we blow up the equilibrium $(0,0)$ of system~\eqref{equ:C7bp} in the positive $u_1$-axis
by the transformation $u_1=u_2^{\alpha_4}$ and $v_1=u_2^{\beta_4}v_2$
and then obtain
\begin{align}
\dot u_2={\cal U}_2(u_2,v_2),~~~
\dot v_2={\cal V}_2(u_2,v_2),
\label{equ:8C7bu2v2p}
\end{align}
where a time rescaling is performed and
{\small
\begin{align*}
{\cal U}_2(u_2,v_2)
&:=u_2\{W_2(u_2,v_2)-c_nu_2^{\delta'-\delta_4}v_2^{\ell-n}+O(u^{\delta'-\delta_4+1})\},
\\
{\cal V}_2(u_2,v_2)
&:=\{\alpha_4(\alpha_3\!-\!\beta_3)\!-\!\beta_4\}v_2W_2(u_2,v_2)
\!+\!(\alpha_4\beta_3\!+\!\beta_4)c_nu_2^{\delta'-\delta_4}v_2^{\ell-n+1}
\!+\!O(u^{\delta'-\delta_4+1}),
\\
W_2(u_2,v_2)
&:=R_{\delta_1}(1,v_2)+u_2R_{\delta_1+1}(1,v_2)+\cdots
+u_2^{\delta'-\delta_4}R_{\delta'}(1,v_2).
\end{align*}
}System~\eqref{equ:8C7bu2v2p} has the equilibrium $(0,0)$ and
at most two equilibria on the $v_2$-axis.
Jacobian matrices at the equilibrium $(0,0)$ and
the equilibrium $(0,v_2^*)$ with $v_2^*\ne 0$ (if exists) are given by
\begin{equation}
\left(
\begin{array}{cccc}
1 & 0
\\
\star & m(n-m)
\end{array}
\right)~~~\mbox{and}~~~
\left(
\begin{array}{cccc}
0 & 0
\\
\star & m(n-m)(\ell-m)b_m(v_2^*)^{\ell-m}
\end{array}
\right),
\label{equ:JC7b1}
\end{equation}
respectively.
We need to determine the property of the semi-hyperbolic equilibrium $(0,v_2^*)$.
By Theorem~7.1 in \cite[p.114]{ZZF},
we solve from the equation ${\cal V}_2(u,v)=0$
near the semi-hyperbolic equilibrium $(0,v_2^*)$
that $v_2=\Lambda_2(u_2)=v_2^*+O(u_2)$ and
{\small
$$
W_2(u_2,\Lambda_2(u_2))=
-\frac{\alpha_4\beta_3+\beta_4}{\alpha_4(\alpha_3-\beta_3)-\beta_4}
c_n(v_2^*)^{\ell-n}u_2^{\delta'-\delta_4}
+O(u_2^{\delta'-\delta_4+1}).
$$
}Substituting it in ${\cal U}_2(u_2,v_2)$ of \eqref{equ:8C7bu2v2p}, we obtain that
$$
{\cal U}_2(u_2,\Lambda_2(u_2))
=-\frac{\alpha_4\alpha_3}{\alpha_4(\alpha_3-\beta_3)-\beta_4}c_n(v_2^*)^{\ell-n} u_2^{\delta'-\delta_4+1}+O(u_2^{\delta'-\delta_4+2}).
$$
Then the qualitative property of the semi-hyperbolic equilibrium $(0,v_2^*)$ (if exists)
of system~\eqref{equ:8C7bu2v2p} is determined by
the parity of $\delta'-\delta_4+1$ and
the signs of
$$
m(n-m)(\ell-m)b_m(v_2^*)^{\ell-m}~~~\mbox{and}~~~
-\frac{\alpha_4\alpha_3}{\alpha_4(\alpha_3-\beta_3)-\beta_4}c_n(v_2^*)^{\ell-n}.
$$
Moreover,
blowing up the equilibrium $(0,0)$ of system~\eqref{equ:C7bp} in the positive $v_1$-axis
by the transformation $u_1=w_2z_2^{\alpha_4}$ and $v_1=z_2^{\beta_4}$,
we obtain
\begin{equation}
\left\{
\begin{array}{llll}
\dot w_2
=\sum_{k\ge \delta_4}\{\beta_4P_k^{(4)}(w_2,1)-\alpha_4w_4Q_k^{(4)}(w_2,1)\}
z_2^{k-\delta_4},
\\
\dot z_2
=z_2\{\sum_{k\ge \delta_4} Q_k^{(4)}(w_2,1)z_2^{k-\delta_4}\},
\end{array}
\right.
\label{equ:8C7bw2z2p}
\end{equation}
When $\alpha_3=\beta_3$, i.e., $m=n+1$,
by {\bf Fact~1},
system~\eqref{equ:8C7bw2z2p} becomes
\begin{equation*}
\dot w_2=w_2\{\beta_4b_m+O(z_2)\},~~~
\dot z_2=z_2\{\beta_3c_nz_2^{\delta'-\delta_4}+O(z_2^{\delta'-\delta_4+1})\}.
\end{equation*}
Then qualitative property of
the semi-hyperbolic equilibrium $(0,0)$ of system~\eqref{equ:8C7bw2z2p}
is determined by the parity of $\delta'-\delta_4$ and the signs of
$\beta_4b_m$ and $\beta_3c_n$.
When $\alpha_3>\beta_3$, i.e., $m>n+1$,
by {\bf Fact~1},
system~\eqref{equ:8C7bw2z2p} becomes
\begin{equation*}
\dot w_2=w_2\{(\beta_4-\alpha_4(\alpha_3-\beta_3))b_m+O(z_2)\},~~~
\dot u_2=z_2\{(\alpha_3-\beta_3)b_m+O(z_2)\}
\end{equation*}
and then the equilibrium $(0,0)$ of system~\eqref{equ:8C7bw2z2p} is
an unstable (or stable) node if $b_m>0$ (or $b_m<0$).

Blowing up the equilibrium $(0,0)$ of system~\eqref{equ:C7bp} in the negative $v_1$-axis
by the transformation $u_1=w_2z_2^{\alpha_4}$ and $v_1=-z_2^{\beta_4}$,
we obtain
\begin{equation}
\left\{
\begin{array}{llll}
\dot w_2
=\sum_{k\ge \delta_4}
\{\beta_4P_k^{(4)}(w_2,-1)+\alpha_4w_2Q_k^{(4)}(w_2,-1)\}z_2^{k-\delta_4},
\\
\dot z_2=-z_2\{\sum_{k\ge \delta_4}Q_k^{(4)}(w_2,-1)z_2^{k-\delta_4}\}.
\end{array}
\right.
\label{equ:8C7bw2z2n}
\end{equation}
When $\alpha_3=\beta_3$, i.e., $m=n+1$,
by {\bf Fact~1},
system~\eqref{equ:8C7bw2z2n} becomes
\begin{equation*}
\dot w_2=w_2\{(-1)^{\ell-m}\beta_4b_m+O(z_2)\},~~~
\dot z_2=z_2\{(-1)^{\delta'-\delta_4}\beta_3c_nz_2^{\delta'-\delta_4}
+O(z_2^{\delta'-\delta_4+1})\}.
\end{equation*}
Then qualitative property of
the semi-hyperbolic equilibrium $(0,0)$ of system~\eqref{equ:8C7bw2z2n}
is determined by the parity of $\delta'-\delta_4$ and the signs of
$(-1)^{\ell-m}\beta_4b_m$ and $(-1)^{\delta'-\delta_4}\beta_3c_n$.
When $\alpha_3>\beta_3$, i.e., $m>n+1$,
by {\bf Fact~1},
system~\eqref{equ:8C7bw2z2n} becomes
\begin{equation*}
\left\{
\begin{array}{llll}
\dot w_2=w_2\{(-1)^{\ell-m}(\beta_4-\alpha_4(\alpha_3-\beta_3))b_m+O(z_2)\},
\\
\dot u_2=z_2\{(-1)^{\ell-m}(\alpha_3-\beta_3)b_m+O(z_2)\}
\end{array}
\right.
\end{equation*}
and its equilibrium $(0,0)$ is a stable (or unstable) node
if $(-1)^{\ell-m}b_m<0$ (or $>0$).

On the other hand,
blowing up the equilibrium $(0,0)$ of system~\eqref{equ:8C1C7t}
in the negative $u$-direction
by the transformation $u=-u_1^{\alpha_3}$ and $v=u_1^{\beta_3}v_1$,
we obtain
\begin{equation}
\left\{
\begin{array}{llll}
\dot u_1=-u_1\frac{1}{(u_1)^{\alpha_3+\delta_3}}\dot u
&=-u_1\{P_{\delta_3}^{(3)}(-1,v_1)+O(u_1)\},
\\
\dot v_1=\frac{\alpha_3}{(u_1)^{\beta_3+\delta_3}}\dot v
+\frac{\beta v_1}{(u_1)^{\alpha_3+\delta_3}}\dot u
&=\widetilde{G}_2(v_1)+O(u_1),
\end{array}
\right.
\label{equ:8C7bu1v1n}
\end{equation}
where $\widetilde{G}_2(v_1):=-v_1^{\ell-m+1}\{\beta_3c_nv_1^{m-n}+(\beta_3-\alpha_3)b_m\}$.
On the $v_1$-axis, system~\eqref{equ:8C7bu1v1n} has the equilibrium $(0,0)$,
which is non-elementary,
and at most two other equilibria.
Jacobian matrix of system~\eqref{equ:8C7bu1v1n} at the equilibrium $(0,\tilde{v}_1^*)$
with $\tilde{v}_1^*\ne 0$ (if exists) is given by
\begin{equation*}
\left(
\begin{array}{cccc}
\frac{\alpha_3}{\beta_3}b_m(\tilde{v}_1^*)^{\ell-m} & 0
\\
\star & \widetilde{G}'_2(\tilde{v}_1^*)
\end{array}
\right).
\end{equation*}
We need further to determine the qualitative property of
the non-elementary equilibrium $(0,0)$ of system~\eqref{equ:8C7bu1v1n}.
Similarly to {\bf Fact~1}, we have the following fact.

\noindent
{\bf Fact~2.}
{\it
The Newton polygon of system~\eqref{equ:8C7bu1v1n} has only one edge,
linking vertex $(0,\alpha_4)$ with vertex $(\beta_4,0)$
and lying on the line $\alpha_4 u+\beta_4 v=\delta_4$.
Moreover,
system~\eqref{equ:8C7bu1v1n} can be rewritten
in quasi-homogeneous components of type $(\alpha_4,\beta_4)$ as
\begin{equation*}
\dot u_1=\sum_{k=\delta_4}^{\delta'}\widetilde{P}_k^{(4)}(u_1,v_1)+\cdots,~~~
\dot v_1=\sum_{k=\delta_4}^{\delta'}\widetilde{Q}_k^{(4)}(u_1,v_1)+\cdots,
\end{equation*}
where dots represent those terms of quasi-homogeneous degree bigger than $\delta'$,
\begin{eqnarray*}
\widetilde{P}_k^{(4)}(u_1,v_1)=u_1\widetilde{R}_k(u_1,v_1),~~~
\widetilde{Q}_k^{(4)}(u_1,v_1)=(\alpha_3-\beta_3)v_1\widetilde{R}_k(u_1,v_1),
\end{eqnarray*}
for all $k=\delta_4,...,\delta'-1$,
\begin{eqnarray*}
&&\widetilde{P}_{\delta'}^{(4)}(u_1,v_1)
=u_1\{c_nv_1^{\ell-n}+\widetilde{R}_{\delta'}(u_1,v_1)\},
\\
&&\widetilde{Q}_{\delta'}^{(4)}(u_1,v_1)
=v_1\{-\beta_3c_nv_1^{\ell-n}
+(\alpha_3-\beta_3)\widetilde{R}_{\delta'}(u_1,v_1)\},
\end{eqnarray*}
and
$
\widetilde{R}_k(u_1,v_1):=b_{\ell-j'_k}u_1^{i'_k}v_1^{j'_k}
+(-1)^{\ell-j''_k+1}a_{\ell-j''_k}u_1^{i''_k}v_1^{j''_k}
$
for all $k=\delta_4,...,\delta'$.
Especially,
$\widetilde{R}_{\delta_4}(u_1,v_1):=b_mv_1^{\ell-m}+(-1)^\ell u_1^{\alpha_3\ell-\delta_3}$.
}

Having {\bf Fact~2},
we blow up the equilibrium $(0,0)$ of system~\eqref{equ:8C7bu1v1n}
in the positive $u_1$-axis
by the transformation $u_1=u_2^{\alpha_4}$ and $v_1=u_2^{\beta_4}v_2$
and then obtain
{\small
\begin{eqnarray}
\left\{
\begin{array}{lllll}
\dot u_2=u_2\{\widetilde{W}_2(u_2,v_2)+c_nu_2^{\delta'-\delta_4}v_2^{\ell-n}
+O(u^{\delta'-\delta_4+1})\},
\\
\dot v_2=\{\alpha_4(\alpha_3-\beta_3)-\beta_4\}v_2\widetilde{W}_2(u_2,v_2)
-(\alpha_4\beta_3+\beta_4)c_nu_2^{\delta'-\delta_4}v_2^{\ell-n+1}
+O(u^{\delta'-\delta_4+1}),
\end{array}
\right.
\label{equ:8C7bu2v2n}
\end{eqnarray}
}where
$
\widetilde{W}_2(u_2,v_2):=\widetilde{R}_{\delta_1}(1,v_2)+u_2\widetilde{R}_{\delta_1+1}(1,v_2)+\cdots
+u_2^{\delta'-\delta_4}\widetilde{R}_{\delta'}(1,v_2).
$
System~\eqref{equ:8C7bu2v2n} has the equilibrium $(0,0)$ and
at most two equilibria on the $v_2$-axis.
Jacobian matrices at the equilibrium $(0,0)$ and
the equilibrium $(0,\tilde{v}_2^*)$ (if exists) with $\tilde{v}_2^*\ne 0$ are given by
\begin{equation}
\left(
\begin{array}{cccc}
(-1)^\ell & 0
\\
\star     & (-1)^{\ell+1}m(m-n)
\end{array}
\right)~\mbox{and}~
\left(
\begin{array}{cccc}
0     & 0
\\
\star & m(n-m)(\ell-m)b_m(\tilde{v}_2^*)^{\ell-m}
\end{array}
\right),
\label{equ:JC7b22}
\end{equation}
respectively.
Similarly,
the qualitative property of the semi-hyperbolic equilibrium $(0,\tilde{v}_2^*)$ (if exists)
of system~\eqref{equ:8C7bu2v2n} is determined by
the parity of $\delta'-\delta_4+1$ and the signs of
$$
m(n-m)(\ell-m)b_m(\tilde{v}_2^*)^{\ell-m}~~~\mbox{and}~~~
\frac{\alpha_4\alpha_3}{\alpha_4(\alpha_3-\beta_3)-\beta_3}c_n(\tilde{v}_2^*)^{\ell-n}.
$$

Moreover,
blowing up the equilibrium $(0,0)$ of system~\eqref{equ:8C7bu1v1n}
in the positive $v_1$-axis
by the transformation $u_1=w_2z_2^{\alpha_4}$ and $v_1=z_2^{\beta_4}$,
we obtain
\begin{equation}
\left\{
\begin{array}{llll}
\dot w_2=\sum_{k\ge\delta_4}
\{\beta_4\widetilde{P}_k^{(4)}(w_2,1)-\alpha_4w_2\widetilde{Q}_k^{(4)}(w_2,1)\}z_2^{k-\delta_4},
\\
\dot z_2=z_2\{\sum_{k\ge\delta_4} \widetilde{Q}_k^{(4)}(w_2,1)z_2^{k-\delta_4}\}.
\end{array}
\right.
\label{equ:8C7bnw2z2p}
\end{equation}
When $\alpha_3=\beta_3$, i.e., $m=n+1$, system~\eqref{equ:8C7bnw2z2p} becomes
\begin{equation*}
\dot w_2=w_2\{\beta_4b_m+O(z_2)\},~~~
\dot z_2=z_2\{-\beta_3c_nz_2^{\delta'-\delta_4}+O(z_2^{\delta'-\delta_4+1})\}.
\end{equation*}
Then qualitative property of the semi-hyperbolic equilibrium $(0,0)$ of system~\eqref{equ:8C7bnw2z2p}
is determined by the parity of $\delta'-\delta_4$ and
the signs of $\beta_4b_m$ and $-\beta_3c_n$.
When $\alpha_3>\beta_3$, i.e., $m>n+1$,
system~\eqref{equ:8C7bnw2z2p} becomes
\begin{equation*}
\dot w_2=w_2\{(\beta_4-\alpha_4(\alpha_3-\beta_3))b_m+O(z_2)\},~~~
\dot u_2=z_2\{(\alpha_3-\beta_3)b_m+O(z_2)\},
\end{equation*}
and its equilibrium $(0,0)$ is an unstable (or stable) node if $b_m>0$ (or $b_m<0$).
Moreover,
blowing up the equilibrium $(0,0)$ of system~\eqref{equ:8C7bu1v1n}
in the negative $v_1$-axis
by the transformation $u_1=w_2z_2^{\alpha_4}$ and $v_1=-z_2^{\beta_4}$,
we obtain
\begin{equation}
\left\{
\begin{array}{llll}
\dot w_2=\sum_{k\ge \delta_4}
\{\beta_4\widetilde{P}_k^{(4)}(w_2,-1)+\alpha_4w_2\widetilde{Q}_k^{(4)}(w_2,-1)\}
z_2^{k-\delta_4},
\\
\dot z_2=-z_2\{\sum_{k\ge \delta_4}\widetilde{Q}_k^{(4)}(w_2,-1)+O(z_2)\}.
\end{array}
\right.
\label{equ:8C7bnw2z2n}
\end{equation}
When $\alpha_3=\beta_3$, i.e., $m=n+1$, system~\eqref{equ:8C7bnw2z2n} becomes
\begin{equation*}
\dot w_2=w_2\{\beta_4b_m+O(z_2)\},~~~
\dot z_2=z_2\{(-1)^{\delta'-\delta_4+1}\beta_3c_nz_2^{\delta'-\delta_4}
+O(z_2^{\delta'-\delta_4+1})\}.
\end{equation*}
Then qualitative property of the semi-hyperbolic equilibrium $(0,0)$ of system~\eqref{equ:8C7bnw2z2n}
is determined by the parity of $\delta'-\delta_4$ and
the signs of $\beta_4b_m$ and $(-1)^{\delta'-\delta_4+1}\beta_3c_n$.
When $\alpha_3>\beta_3$, i.e., $m>n+1$,
system~\eqref{equ:8C7bnw2z2n} becomes
\begin{equation*}
\left\{
\begin{array}{llll}
\dot w_2=w_2\{(-1)^{\ell-m}(\beta_4-\alpha_4(\alpha_3-\beta_3))b_m+O(z_2)\},
\\
\dot u_2=z_2\{(-1)^{\ell-m}(\alpha_3-\beta_3)b_m+O(z_2)\},
\end{array}
\right.
\end{equation*}
and its equilibrium $(0,0)$ is an unstable (or stable) node
if $(-1)^{\ell-m}b_m>0$ (or $(-1)^{\ell-m}b_m<0$).

\begin{figure}[H]
\centering
\subcaptionbox{%
     }{\includegraphics[height=1in]{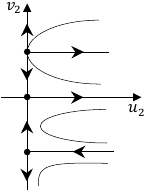}}~
\subcaptionbox{%
     }{\includegraphics[height=1in]{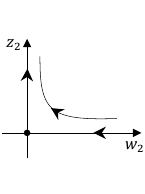}}~
\subcaptionbox{%
     }{\includegraphics[height=1in]{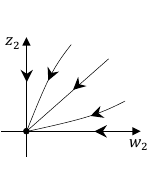}}~
\subcaptionbox{%
     }{\includegraphics[height=1in]{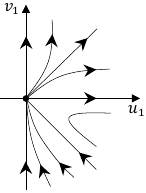}}
\\
\subcaptionbox{%
     }{\includegraphics[height=1in]{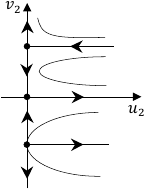}}~
\subcaptionbox{%
     }{\includegraphics[height=1in]{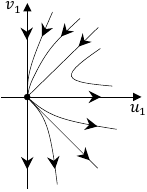}}~
\subcaptionbox{%
     }{\includegraphics[height=1in]{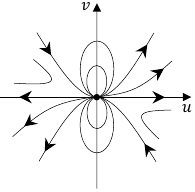}}
\caption{Phase portraits in {\bf(C7b)} with $m=n+1$.}
\label{fig:C7b1}
\end{figure}

When $m=n+1$,
only situations {\bf(S3)}-{\bf(S5)} and {\bf(S8)}-{\bf(S10)} need to be considered.
We only give the proof in {\bf(S10)} since proofs in other situations are similar.
In {\bf(S10)},
system~\eqref{equ:8C7bu2v2p} has three equilibria on the $v_2$-axis and
their properties can be determined by \eqref{equ:JC7b1}.
The phase portrait of system~\eqref{equ:8C7bu2v2p} along the $v_2$-axis is
given by Fig.~\ref{fig:C7b1}(a).
Phase portraits of systems~\eqref{equ:8C7bw2z2p} and \eqref{equ:8C7bw2z2n}
near the origin are given by Fig.~\ref{fig:C7b1}(b) and (c), respectively.
Note that $(0,0)$ is the only equilibrium of system~\eqref{equ:C7bp} on the $v_1$-axis.
So, after blowing down,
we obtain the phase portrait Fig.~\ref{fig:C7b1}(d) of system~\eqref{equ:C7bp}
along the $v_1$-axis.
On the other hand,
system~\eqref{equ:8C7bu2v2n} has three equilibria on the $v_2$-axis and
their properties can be determined by \eqref{equ:JC7b22}.
So we obtain the phase portrait Fig.~\ref{fig:C7b1}(e) of system~\eqref{equ:8C7bu2v2n}
along the $v_2$-axis.
Moreover,
phase portraits of systems~\eqref{equ:8C7bnw2z2p} and \eqref{equ:8C7bnw2z2n}
near the origin are
given by Fig.~\ref{fig:C7b1}(c) and (b), respectively.
Note that system~\eqref{equ:8C7bu1v1n} has the only equilibrium $(0,0)$ on the $v_1$-axis.
So, after blowing down,
we obtain the phase portrait Fig.~\ref{fig:C7b1}(f) of system~\eqref{equ:8C7bu1v1n}
along the $v_1$-axis.
Combining Fig~\ref{fig:C7b1} (d) and (f) and blowing down,
we obtain the phase portrait Fig.~\ref{fig:C7b1}(g) of system~\eqref{equ:8C1C7t}.
So the phase portrait of system~\eqref{GL:pfg} near the equator of the Poincar\'{e} disc
is given by Fig.~\ref{fig:infty}(s).

\begin{figure}[H]
\centering
\subcaptionbox{%
     }{\includegraphics[height=1in]{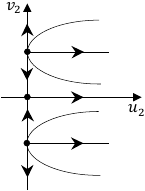}}~
\subcaptionbox{%
     }{\includegraphics[height=1in]{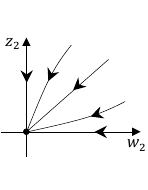}}~
\subcaptionbox{%
     }{\includegraphics[height=1in]{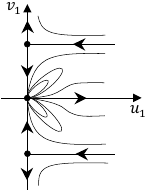}}~
\subcaptionbox{%
     }{\includegraphics[height=1in]{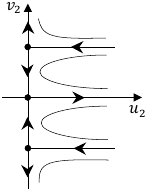}}~
\subcaptionbox{%
     }{\includegraphics[height=1in]{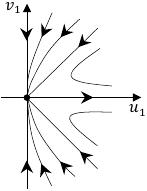}}~
\subcaptionbox{%
     }{\includegraphics[height=1in]{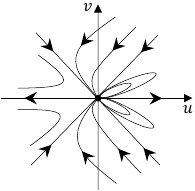}}
\caption{Phase portraits in {\bf(C7b)} with $m>n+1$.}
\label{fig:C7b2}
\end{figure}

When $m>n+1$,
there are 12 situations {\bf(S1)}-{\bf(S12)} and we only give the proof in {\bf(S12)}
since proofs in other situations are similar.
In {\bf(S12)},
system~\eqref{equ:8C7bu2v2p} has three equilibria on the $v_2$-axis and
the phase portrait along the $v_2$-axis is given by Fig.~\ref{fig:C7b2}(a).
Moreover,
phase portrait of systems~\eqref{equ:8C7bw2z2p} and \eqref{equ:8C7bw2z2n} are
both given by Fig.~\ref{fig:C7b2}(b).
System~\eqref{equ:C7bp} has three equilibria on the $v_1$-axis,
and properties of equilibria not at the origin can be determined by \eqref{equ:JC7b21}.
After blowing down,
we obtain from phase portraits of systems~\eqref{equ:8C7bu2v2p},
\eqref{equ:8C7bw2z2p} and \eqref{equ:8C7bw2z2n}
the phase portrait Fig.~\ref{fig:C7b2}(c) of system~\eqref{equ:C7bp} along the $v_1$-axis.
On the other hand,
phase portrait of system~\eqref{equ:8C7bu2v2n} along the $v_2$-axis is
given by Fig.~\ref{fig:C7b2}(d),
and phase portraits of systems~\eqref{equ:8C7bnw2z2p} and \eqref{equ:8C7bnw2z2n} are both
given by Fig.~\ref{fig:C7b2}(b).
So after blowing down,
we obtain phase portrait Fig.~\ref{fig:C7b2}(e) of system~\eqref{equ:8C7bu1v1n}
along the $v_1$-axis.
Combining Fig.~\ref{fig:C7b2}(c) and (e) and blowing down,
we obtain the phase portrait Fig.~\ref{fig:C7b2}(f) of system~\eqref{equ:8C1C7t}
near the origin.
Thus the phase portrait of system~\eqref{GL:pfg} near the equator of the Poincar\'{e} disc
is given by Fig.~\ref{fig:infty}(s).
Consequently, we complete the proof of Theorem~\ref{thm:infty}.
\qquad$\Box$

\section{Global center and its non-isochronicity}
\setcounter{equation}{0}
\setcounter{lm}{0}
\setcounter{thm}{0}
\setcounter{rmk}{0}
\setcounter{df}{0}
\setcounter{cor}{0}
\setcounter{pro}{0}

As an application of the above two sections,
we classify all global phase portraits of system~\eqref{GL:pfg}
when the origin is the only equilibrium and a center.
We further prove the non-isochronicity of the global center
by investigating period of orbit near the equator.

\begin{figure}[H]
\centering
\subcaptionbox{%
     }{\includegraphics[height=1in]{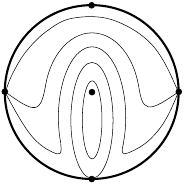}}
\subcaptionbox{%
     }{\includegraphics[height=1in]{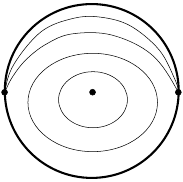}}
\subcaptionbox{%
     }{\includegraphics[height=1in]{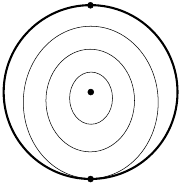}}
\subcaptionbox{%
     }{\includegraphics[height=1in]{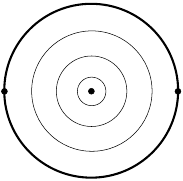}}
\subcaptionbox{%
     }{\includegraphics[height=1in]{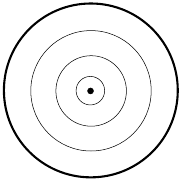}}
\caption{Global phase portraits of system~\eqref{GL:pfg} with $\ell\ge 2$
when the origin $O$ is the only finite equilibrium and a center.}
\label{fig:OG}
\end{figure}

\begin{thm}
{\bf (i)}
If the origin is the only equilibrium and a center of system~\eqref{GL:pfg} with $\ell\ge 2$,
then there are 5 different topological global phase portraits,
see Fig.~\ref{fig:OG}.
\\
{\bf (ii)}
System~\eqref{GL:pfg} with $\ell\ge 2$ and $a_p=a_\ell=-1$ has a global center at the origin,
i.e., having global phase portrait Fig.~\ref{fig:OG} (d) or (e),
if and only if
\begin{description}
  \item[(G1)]
neither $\varphi$ nor $g$ has nonzero real roots,

  \item[(G2)]
condition {\bf(M1)} or {\bf(M2)}, given in Theorem~\ref{thm:center}, holds,

  \item[(G3)]
the system $F(x)=F(z)$ and $G(x)=G(z)$
has a unique solution $z(x)$ satisfying $z(0)=0$ and $z'(0)<0$,
where $G(x):=\int_0^x g(s) ds$,

  \item[(G4)]
condition {\bf(W1)}, {\bf(W2)} or {\bf(W3)}, given in Remark~\ref{rmk:MI}, holds.
\end{description}
{\bf(iii)} Global center of system~\eqref{GL:pfg} cannot be isochronous.
\label{thm:global}
\end{thm}

\begin{rmk}
Due to the transformation $y\to \alpha y$ and $t\to \frac{-1}{a_p\alpha^p}t$
with $\alpha:=\big|\frac{a_p}{a_\ell}\big|^{\frac{1}{\ell-p}}$,
we can assume without loss of generality that
$(a_p,a_\ell)=(-1,\pm1)$ in system~\eqref{GL:pfg}.
When $(a_p,a_\ell)=(-1,1)$, the polynomial $\varphi$ has a nonzero real root,
implying that the origin is not the only equilibrium and therefore
there is no global center.
\end{rmk}

\begin{eg}
Global phase portraits Fig.~\ref{fig:OG} (a)-(e) can be realized by
the following generalized polynomial Li\'{e}nard systems
\begin{equation*}
\begin{array}{lllll}
&\dot x=-y^3-x^4,     &\dot y=x^3,
\\
&\dot x=-y^5-x^4,     &\dot y=x^3,
\\
&\dot x=-y^3-x^{10},  &\dot y=x^{11},
\\
&\dot x=-y^3-y^7-x^4, &\dot y=x^3+x^5,
\\
&\dot x=-y^3-y^5-x^4, &\dot y=x^3+x^5,
\end{array}
\end{equation*}
respectively.
\end{eg}

{\bf Proof of Theorem~\ref{thm:global}.}
For {\bf(i)},
since the origin $O$ is the only equilibrium and a center,
its period annulus ${\cal U}$ is unbounded and therefore
either ${\cal U}=\mathbb{R}^2\setminus\{O\}$, or
${\cal U}\varsubsetneqq\mathbb{R}^2\setminus\{O\}$.
In the first case,
we obtain global phase portraits Fig.~\ref{fig:OG} (d) and (e)
from Fig.~\ref{fig:infty} (w) and (x).

In the second case,
there is an infinite equilibrium having a hyperbolic sector
whose two separatrices are contained in the boundary of the period annulus ${\cal U}$
and moreover at least one separatrix does not lie on the equator.
This implies directly that
the phase portrait near infinity can not be
Fig.~\ref{fig:infty} (a), (d), (i), (k), (n), (o), (w) and (x).

Moreover,
we claim that
the phase portrait near infinity can not be
Fig.~\ref{fig:infty} (b), (c), (e), (f), (g), (h), (p), (r) and (s).
In fact, in each one of those phase portraits,
one separatrix of each hyperbolic sector links with a parabolic sector,
implying that this separatrix can not be the boundary of the period annulus ${\cal U}$.

Finally,
we claim that
the phase portrait near infinity can not be
Fig.~\ref{fig:infty} (j), (m), (t) and (v).
Actually,
the Poincar\'{e}-Hopf Index Theorem (\cite[Theorem~6.30]{DLA}) on the sphere $\mathbb{S}^2$ states that for every tangent vector field on $\mathbb{S}^2$ with a finite number of equilibria,
the sum of their indices is 2.
We use the Poincar\'{e} sphere compactification to
extend the vector field generated by system~\eqref{GL:pfg} to
a vector field on the sphere $\mathbb{S}^2$.
In the case of Fig.~\ref{fig:infty} (j),
there are totally six equilibria on the the sphere $\mathbb{S}^2$,
two centers at the south and the north and
two nodes and four saddle-nodes on the equator.
So the sum of indices of all equilibria are 4,
a contradiction.
It is similar to analyze Fig.~\ref{fig:infty} (m), (t) and (v).

Consequently,
in the second case ${\cal U}\varsubsetneqq\mathbb{R}^2\setminus\{O\}$,
the phase portrait near infinity can only be
one of Fig.~\ref{fig:infty} (l), (q) and (u).
So there are only 5 different topological phase portraits and
Theorem~\ref{thm:global}{\bf(i)} is proved.

For {\bf(ii)},
as indicated in \cite[Proposition 2]{LV21},
a polynomial differential system has a global center at the origin
if and only if the origin is the only equilibrium and is a local center and
the system is monodromic at infinity.
For the sufficiency,
condition {\bf(G1)} clearly implies that the origin is the only equilibrium,
conditions {\bf(G2)} and {\bf(G3)} ensure that
the origin is a center by Theorem~\ref{thm:center},
and
condition {\bf(G4)} ensures that
system~\eqref{GL:pfg} is monodromic at infinity because of Remark~\ref{rmk:MI}.
For the necessity,
we see from Proposition 2 in \cite{LV21}, Theorem~\ref{thm:center} and Remark~\ref{rmk:MI}
that conditions {\bf(G2)}, {\bf(G3)} and {\bf(G4)} holds.
If condition {\bf(G1)} is invalid,
then either $\varphi$ has a nonzero real root $y_*$,
or $g$ has a nonzero real root $x_*$.
In the first case, $(0,y_*)$ is an equilibrium of system~\eqref{GL:pfg}, a contradiction.
In the second case,
the polynomial $\varphi(y)-F(x_*)$ has a nonzero real root $\hat y$
because its degree $\ell$ is odd as required in {\bf(G4)},
implying that $(x_*,\hat y)$ is an equilibrium of system~\eqref{GL:pfg}, a contradiction.
Hence, condition {\bf(G1)} holds and Theorem~\ref{thm:global}{\bf (ii)} is proved.

For {\bf(iii)},
global center only appears in cases {\bf(C1)}, {\bf(C3)}, {\bf(C4)} and {\bf(C7)}
by Theorem~\ref{thm:infty}.
Moreover,
phase portrait of global center
in cases {\bf(C1)} and {\bf(C3)} is given by Fig.~\ref{fig:OG} (e),
and that in cases {\bf(C4)} and {\bf(C7)} is given by Fig.~\ref{fig:OG} (d).
As indicated in Theorem~2.6 of \cite{CD97},
if a nonlinear differential system has an isochronous center,
then the system has at least one equilibrium on the equator of the Poincar\'{e} disc.
So a global center in cases {\bf(C1)} and {\bf(C3)} is clearly not isochronous.

\begin{figure}[H]
\centering
\subcaptionbox{%
     }{\includegraphics[height=1in]{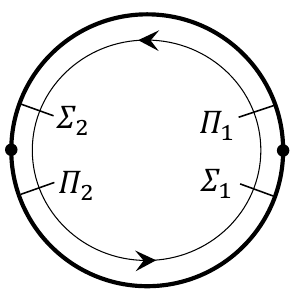}}~~~~~~
\subcaptionbox{%
     }{\includegraphics[height=1in]{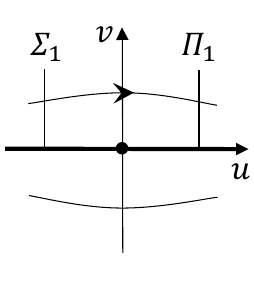}}~~~~~~
\subcaptionbox{%
     }{\includegraphics[height=1in]{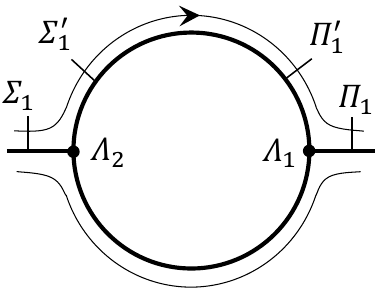}}
\caption{(a) Transverse sections near equilibria of
             system~\eqref{GL:pfg} on the equator.
         (b) Transverse sections near equilibrium $\theta=0$ in the local chart $(u,v)$.
         (c) Desingularization of the equilibrium $\theta=0$ and
             transverse sections near the divisor.}
\label{fig:2double}
\end{figure}

In what follows,
we only prove the non-isochronicity of the global center in case {\bf(C7)}
since the proof in case {\bf(C4)} is similar.
In this case,
system~\eqref{GL:pfg} has only two equilibria
$\theta_1:=0$ and $\theta_2:=\pi$ on the equator.
Let $\zeta_k,\xi_k:[0,1]\to\mathbb{D}^2$ be analytic curves transverse to
the equator near the infinite equilibrium $\theta_k$
such that $\zeta_k(0)<\theta_k<\xi_k(0)$ on the equator.
Let $\Sigma_k$ and $\Pi_k$ denote the images of $\zeta_k$ and $\xi_k$, respectively.
Then the equator is divided into two corners at $\theta_k$
and two sides $S_k$ between the corners at $\theta_{k}$ and $\theta_{k+1}$
with $\theta_3:=\theta_1$, as shown in Fig.~\ref{fig:2double} (a).

Suppose that $\varphi(t,\zeta_1(s))$ is the orbit
passing through the point $\zeta_1(s)\in\Sigma_1$ for small $s>0$.
Let $t_1(s)$ be the passage time near the corner at $\theta_1$ on the equator,
i.e.,
$$
t_1(s):=\min\{t\in\mathbb{R}_+:\varphi(t,\zeta_1(s))\in \Pi_1\}.
$$
Similarly,
we can define the passage time $t_2(s)$ near the corner at $\theta_2$ on the equator,
and the passage time $\tau_k(s)$ along the side $S_k$ for $k=1$ and $2$.
We claim that
\begin{align}
\lim_{s\to 0^+}t_k(s)=\lim_{s\to 0^+}\tau_k(s)=0~~~\forall k=1,2.
\label{tstaos}
\end{align}
Actually,
we only need to consider the side $S_1$ and the corner at $\theta_1$
since it is similar to investigate the others.
Along the side $S_1$,
the vector field associated to system~\eqref{GL:pfg} is of the form
$$
\widetilde{\cal X}^{(0)}
:=\frac{1}{v^{\ell-1}}{\cal X}^{(0)}
=\frac{1}{v^{\ell-1}}
\left\{
\big(-1+O(v)\big)\frac{\partial}{\partial u}
+O(v)\frac{\partial}{\partial v}
\right\}
$$
for $u\in(-M,M)$,
where ${\cal X}^{(0)}$ is the vector field generated by system~\eqref{equ:8C2C7t},
$M$ is a large positive constant.
Since $\ell-1>0$,
Lemma~2.4 of \cite{LSZ} ensures that
$$
\lim_{s\rightarrow 0^+}\tau_1(s)=0.
$$
For the corner at $\theta_1$,
the vector field associated to system~\eqref{GL:pfg} is of the form
$$
\widetilde{\cal Y}^{(0)}
:=\frac{1}{v^{\ell-1}}{\cal Y}^{(0)},
$$
where ${\cal Y}^{(0)}$ is the vector field generated by system~\eqref{equ:8C1C7t}.
After desingularization,
the hyperbolic sector of the equilibrium $(0,0)$ of $\widetilde{\cal Y}^{(0)}$
in the half-plane $v>0$ becomes
a poly-arc with two hyperbolic saddles $\Lambda_1$ and $\Lambda_2$ as vertices,
as shown in Fig.~\ref{fig:2double}(c).
Similarly to the above,
let $\Gamma_1$ denote the side lying between the two corners at $\Lambda_1$ and $\Lambda_2$.
In the blow-up coordinates,
the desingularized vector field of $\widetilde{\cal Y}^{(0)}$
near the corner at $\Lambda_1$ takes the form
$$
\frac{u_1^{\delta_2}{\cal Y}_1^{(1)}}{\alpha_2(u_1^{\beta_2}v_1)^{\ell-1}}
=\frac{1}{\alpha_2u_1^{\ell n-1}v_1^{\ell-1}}
\left\{
u_1(1+o(1))\frac{\partial }{\partial u_1}
-v_1\left(n+1+o(1)\right)\frac{\partial }{\partial v_1}
\right\},
$$
where ${\cal Y}_1^{(1)}$ is the vector field generated by system~\eqref{equ:8C4ap}.
In the blow-up coordinates,
the desingularized vector field of $\widetilde{\cal Y}^{(0)}$
near the corner at $\Lambda_2$ takes the form
$$
\frac{u_1^{\delta_2}{\cal Y}_2^{(1)}}{\alpha_2(u_1^{\beta_2}v_1)^{\ell-1}}
=\frac{(-1)^\ell}{\alpha_2u_1^{\ell n-1}v_1^{\ell-1}}
\left\{
u_1(1+o(1))\frac{\partial }{\partial u_1}
-v_1(n+1+o(1))\frac{\partial }{\partial v_1}
\right\},
$$
where ${\cal Y}_2^{(1)}$ is the vector field generated by system~\eqref{equ:8C4an}.
Since $\ell n-1>0$ and $\ell-1>0$,
by Lemma~2.4 of \cite{LSZ},
the passage times near the two hyperbolic saddle corners at $\Lambda_1$ and $\Lambda_2$
both approach to $0$.
On the other hand,
applying the transformation $u=w_1z_1^{\alpha_2}$ and $v=z_1^{\beta_2}$ to
the vector field $\widetilde{\cal Y}^{(0)}$,
we obtain that
the desingularized vector field of $\widetilde{\cal Y}^{(0)}$
along the side $\Gamma_1$ is of the form
$$
\frac{1}{\beta_2z_1^{\ell n-1}}
\left\{
\big(H(w_1)+O(z_1)\big)\frac{\partial }{\partial w_1}
-z_1\big(w_1^\ell+O(z_1)\big)\frac{\partial }{\partial z_1}
\right\},
$$
where $H(w_1):=(n+1)w_1^{\ell+1}-\beta_2c_n$.
We see from Theorem~\ref{thm:infty} that
if system~\eqref{GL:pfg} has a global center in case {\bf(C7)},
then condition {\bf(T2)} also holds.
Condition {\bf(T2)} implies that $H(w_1)$ has no real roots.
So Lemma~2.4 of \cite{LSZ} ensures that
the passage time along the side $\Gamma_1$ approaches to $0$.
Consequently,
we obtain that
$$
\lim_{s \to 0^+} t_1(s)=0.
$$
Thus our claimed \eqref{tstaos} is proved.
It follows that
if system~\eqref{GL:pfg} has a global center in case {\bf(C7)},
then the period of orbit near the infinity approaches to $0$
and therefore global center is not isochronous.
Thus the proof of Theorem~\ref{thm:global} is completed.
\qquad$\Box$



\end{document}